\documentclass{amsart}
\usepackage{amssymb}
\usepackage{graphicx,color}
\usepackage{amscd}
\usepackage{stmaryrd}

\calclayout
\allowdisplaybreaks

\newcommand{\ignore}[1]{}

\numberwithin{figure}{section}
\numberwithin{table}{section}
\newcommand\tr{\operatorname{tr}}

\newcommand\grad{\operatorname{grad}}
\renewcommand\div{\operatorname{div}}
\newcommand\curl{\operatorname{curl}}
\newcommand\rot{\operatorname{rot}}

\newcommand\spn{\operatorname{span}}

\newcommand\N{\mathbb{N}}
\newcommand\R{\mathbb{R}}

\newcommand\x{\times}

\newcommand\D{d}  

\renewcommand\H{\mathcal H}

\newcommand\I{{\mathcal I}}

\renewcommand\P{{\mathcal P}}

\newcommand\T{{\mathcal T}}

\newcommand\W{{\mathbb W}}

\newcommand\Tr{\operatorname{Tr}}

\newcommand{\0}{\mathaccent23}

\newcommand\Alt{\operatorname{Alt}}

\newcommand\Pfalpha{P_{f,\alpha}}

\numberwithin{equation}{section}

\newtheorem{thm}{Theorem}[section]
\newtheorem{prop}[thm]{Proposition}
\newtheorem{lem}[thm]{Lemma}
\newtheorem{cor}[thm]{Corollary}

\newenvironment{remark}[1]
{\medskip \noindent {\bf Remark.} #1}

\begin{document}

\title[Geometric decompositions]{Geometric decompositions and local bases for spaces of finite element
  differential forms}
%    author one information
%\author[D. N. Arnold, R. S. Falk, and R. Winther]
\author{Douglas N. Arnold}
\address{Institute for Mathematics and its Applications
and School of Mathematics,
University of Minnesota, Minneapolis, MN 55455}
\email{arnold@ima.umn.edu}
\urladdr{http://www.ima.umn.edu/\char'176arnold/}
\thanks{}
%    author two information
\author{Richard S. Falk}
\address{Department of Mathematics,
Rutgers University, Piscataway, NJ 08854}
\email{falk@math.rutgers.edu}
\urladdr{http://www.math.rutgers.edu/\char'176falk/}
\thanks{}
%    author three information
\author{Ragnar Winther}

\address{Centre of Mathematics for Applications
and Department of Informatics,
University of Oslo, 0316 Oslo, Norway}
\email{ragnar.winther@cma.uio.no}
\urladdr{http://heim.ifi.uio.no/\char'176rwinther/}
\thanks{}
%    \subjclass is required.
\date{May 12, 2008}
\thanks{The work of the first author was supported in part by NSF grant
DMS-0713568.  The work of the second author was supported in part by NSF grant
DMS06-09755. The work of the third author was supported by the Norwegian
Research Council.}

\begin{abstract}
We study the two primary families of spaces of finite element
differential forms with respect to a simplicial mesh in any number of space
dimensions. These spaces
are generalizations of the classical finite element spaces for vector fields,
frequently referred to as Raviart--Thomas, Brezzi--Douglas--Marini,
and N\'ed\'elec spaces. In the present paper, we derive
geometric decompositions of these spaces which lead directly to explicit
local bases for them, generalizing the Bernstein
basis for ordinary Lagrange finite elements. The approach applies to both
families of finite element spaces, for arbitrary polynomial degree, arbitrary
order of the differential forms, and an arbitrary simplicial triangulation
in any number of space dimensions.  A prominent role in the construction
is played by the notion of a consistent family of extension operators, which
expresses in an abstract framework a sufficient condition for deriving a
geometric decomposition of a finite element space leading to a local basis.
\end{abstract}

\maketitle

\section{Introduction}\label{sec:intro}

The study of finite element exterior calculus has given increased insight into
the construction of stable and accurate finite element methods for problems
appearing in various applications, ranging from electromagnetics to
elasticity. Instead of considering the design of discrete methods for each
particular problem separately, it has proved beneficial to simultaneously
study approximations of a family of problems, tied together by a common
differential complex.

To be more specific, let $\Omega \subset \R^n$ and 
let $H\Lambda^k(\Omega)$ be the space of differential $k$ forms 
$\omega$ on $\Omega$,
which is in $L^2$, and where its exterior derivative,
$d\omega$,
is also in $L^2$.
% This space is a Hilbert space.
The $L^2$ version of the de Rham complex 
then takes the form
 \[
0\to H\Lambda^0(\Omega)
\xrightarrow{d} H\Lambda^1(\Omega) \xrightarrow{d} 
\cdots \xrightarrow{d}
H\Lambda^n(\Omega)\to0.
\]
The basic construction in finite element exterior calculus is
of a corresponding subcomplex
\[
0\to\Lambda_h^0\xrightarrow{d} \Lambda_h^1 \xrightarrow{d} 
\cdots \xrightarrow{d}
\Lambda_h^n\to0,
\]
where the spaces $\Lambda_h^k$ are finite dimensional subspaces of
$H\Lambda^k(\Omega)$ consisting
of piecewise polynomial differential forms with respect to a partition
of the domain $\Omega$. In the theoretical analysis of the stability of numerical methods
constructed from this discrete complex, bounded projections $\Pi_h :
H\Lambda^k(\Omega) \to \Lambda_h^k$ are utilized, such that the diagram
\begin{equation*}
\begin{CD}
0\to@.H\Lambda^0(\Omega) @>d>> H\Lambda^1(\Omega) @>d>>
\cdots @>d>> H\Lambda^n(\Omega)@.\to0 \\
@. @VV\Pi_hV @VV\Pi_hV @.  @VV\Pi_hV @.\\
0\to@.\Lambda^0_h @>d>> 
\Lambda^1_h @>d>>
\cdots @>d>> \Lambda^n_h @.\to0
\end{CD}
\end{equation*}
commutes. For a general reference to finite element
exterior calculus, we refer to the survey paper \cite{acta}, and references
given therein.  As is shown there, the spaces $\Lambda^k_h$ are
taken from two main families.  Either $\Lambda^k_h$ is one
of the spaces $\P_r\Lambda^k(\T)$ consisting of all elements
of $H\Lambda^k(\Omega)$ which restrict to polynomial $k$-forms
of degree at most $r$ on each simplex $T$ in the partition
$\T$, or $\Lambda^k_h=\P^-_r\Lambda^k(\T)$, which is a space
which sits between $\P_r\Lambda^k(\T)$ and $\P_{r-1}\Lambda^k(\T)$
(the exact definition will be recalled below).
These spaces are
generalizations of the Raviart-Thomas and
Brezzi-Douglas-Marini spaces used to discretize $H(\div)$ and $H(\rot)$ in two space dimensions and the N\'ed\'elec edge
and face spaces of the first and second kind, used to
discretize $H(\curl)$ and $H(\div)$ in three space dimensions.

A key aim of the present paper is to explicitly construct geometric
decompositions of the spaces $\P_r\Lambda^k(\T)$ and
$\P^-_r\Lambda^k(\T)$ for arbitrary values of $r\ge1$ and $k\ge0$, and
an arbitrary simplicial partition $\T$ of a polyhedral domain in an
arbitrary number of space dimensions.  More precisely, we will
decompose the space into a direct sum with summands indexed by the
faces of the mesh (of arbitrary dimension), such that the summand
associated to a face is the image under an explicit extension operator
of a finite-dimensional space of differential forms on the face.  Such
a decomposition is necessary for an efficient implementation of the finite
element method, since it allows an assembly process that leads to local
bases for the finite element space. The construction of explicit local
bases is the other key aim of this work.

The construction given here leads to a generalization of the
so-called Bernstein basis for ordinary polynomials, i.e., 0-forms on a simplex
$T$ in $\R^n$, and the corresponding finite element spaces, the Lagrange finite elements.  See Section~\ref{Bernstein} below.  This polynomial basis is a
well known and useful theoretical tool both in finite element analysis and
computational geometry. For low order piecewise polynomial spaces, it can be
used directly as a computational basis, while for polynomials of higher order,
this basis can be used as a starting point to construct a basis with improved
conditioning or other desired properties.  The same will be true for the
corresponding bases for spaces of piecewise polynomial differential forms
studied in this paper.

This paper continues the development of geometric decompositions
begun in \cite[Section 4]{acta}.  In the present paper, we give a prominent place to
the notion of a \emph{consistent family of extension
operators}, and show that such a family leads to
a direct sum decomposition of the piecewise polynomial space of
differential forms with proper interelement
continuity.  The explicit notion of a consistent family of extension
operators is new to this paper.  We also take a more geometric and coordinate-independent
approach in this paper than in \cite{acta}, and so are able to give a
purely geometric characterization of the decompositions
obtained here. The geometric decomposition
we present for the spaces $\P^-_r \Lambda^k$ here turns out to be the same as
obtained in \cite{acta}, but the decomposition of the spaces $\P_r
\Lambda^k$ obtained here is new.  It improves upon the one obtained in
\cite{acta}, since it no longer depends on a particular choice of
ordering of the vertices of the simplex $T$, and leads to a more canonical for $\P_r \Lambda^k$.

The construction of implementable bases for some of the spaces we
consider here has been considered previously by a number of authors.
Closest to the present paper is the work of Gopalakrishnan,
Garc\'{\i}a-Castillo, and Demkowicz \cite{g-g-d}.  They give a basis
in barycentric coordinates for the space $\P^-_r\Lambda^1$, where
$T$ is a simplex in any number of space dimensions. In this particular
case, their basis is the same as we present in
Section~\ref{sec:bases}. In fact, Table 3.1 of \cite{g-g-d} is the
same, up to a change in notation, as the left portion of Table
\ref{tb:t2} of this paper. As will be seen below, explicit bases for
the complete polynomial spaces $\P_r\Lambda^k$ are more complicated
than for the $\P^-_r\Lambda^k$ spaces. To our knowledge, the basis
we present here for the $\P_r\Lambda^k$ spaces have not previously
appeared in the literature, even in two dimensions or for small values of $r$.

Other authors have focused on the construction of $p$-hierarchical
bases for some of the spaces considered here. We particularly note the
work of Ainsworth-Coyle \cite{ainsworth-coyle}, Hiptmair
\cite{hiptmair-pier}, and Webb \cite{webb}. In \cite{ainsworth-coyle},
the authors construct hierarchical bases of arbitrary polynomial order
for the spaces we denote $\P_r\Lambda^k$, $k=0, \ldots, 3$, $r \ge
1$, and $T$ a simplex in three dimensions.  In section 5 of
\cite{hiptmair-pier}, Hiptmair considers hierarchical bases of
$\P^-_r\Lambda^k$ for general $r$, $k$, and simplex dimension.  In
\cite{webb}, Webb constructs hierarchical bases for both
$\P_r\Lambda^k$ and $\P^-_r\Lambda^k$, for $k=0,1$ in one,
two,and three space dimensions.  The approaches of these three sets of
authors differ. Even when adapted to the simple case of zero-forms,
i.e., Lagrange finite elements, they produce different hierarchical
bases, from among the many that have been proposed.  Our approach is
quite distinct from these in that we are not trying to find
hierarchical bases, but rather we generalize the explicit Bernstein
basis to the full range of spaces $\P_r\Lambda^k$ and
$\P^-_r\Lambda^k$.

In the present work, by treating
the $\P_r\Lambda^k$ and $\P^-_r\Lambda^k$ families
together, and adopting the framework of
differential forms, we are able to give a presentation that shows the close
connection of these two families, and is valid for all order polynomials and
all order differential forms in arbitrary space dimensions.  Moreover, the
viewpoint of this paper is that the construction of basis functions is a
straightforward consequence of the geometric decomposition of the finite
element spaces, which is the key ingredient needed to construct spaces with the
proper inter-element continuity.  Thus, the main results of the paper
focus on these geometric decompositions.

An outline of the paper is as follows.  In the next section, we define
our notation and review material we will need about barycentric
coordinates, the Bernstein basis, differential forms, and simplicial
triangulations.  The $\P_r\Lambda^k$ and $\P^-_r\Lambda^k$ families of
polynomial and piecewise polynomial differential forms are described
in Section~\ref{sec:pdf}. In Section~\ref{cgd}, we introduce the
concept of a consistent family of extension operators and use it to construct
a geometric decomposition of a finite element space in an abstract setting. In addition to the Bernstein
decomposition, a second familiar decomposition which fits this framework is the dual
decomposition, briefly discussed in Section~\ref{sec:dof}. Barycentric spanning sets and bases for the spaces $\P_r \Lambda^k(T)$
and $\P^-_r \Lambda^k(T)$ and the corresponding subspaces $\0\P_r \Lambda^k(T)$
and $\0\P^-_r \Lambda^k(T)$
with vanishing trace are presented in
Section~\ref{sec:span-basis}. The main results of this paper, the
geometric decompositions and local bases, are derived in Sections~\ref{sec:decomp-}
and \ref{sec:decomp} for $\P_r^- \Lambda^k(\T)$
and $\P_r \Lambda^k(\T)$.  Finally, in Section~\ref{sec:bases}, we discuss how
these results can be used to obtain explicit local bases, and tabulate
such bases in the cases of $2$ and $3$ space
dimensions and polynomial degree at most $3$.

\section{Notation and Preliminaries}\label{sec:not-prelim}
\subsection{Increasing sequences and multi-indices}
\label{sec:notation}
We will frequently use increasing sequences, or increasing maps from integers
to integers, to index differential forms.  For integers $j,k,l,m$, with
$0\le k-j \le m-l$, we will use $\Sigma(j:k,l:m)$ to denote
the set of increasing maps $\{j,\ldots,k\}\to\{l,\ldots,m\}$, i.e.,
\begin{equation*}
\Sigma(j:k,l:m) = \{\, \sigma:\{j,\ldots,k\}\to\{l,\ldots,m\}\,|\,
\sigma(j)<\sigma(j+1)< \cdots<\sigma(k)\,\}.
\end{equation*}
Furthermore, $\llbracket \sigma \rrbracket$ will denote the
range of such maps, i.e., for $\sigma\in\Sigma(j:k,l:m)$,
$\llbracket\sigma\rrbracket = \{\,\sigma(i)\,|\, i=j,\ldots,k\}$.
Most frequently, we will use the sets
$\Sigma(0:k,0:n)$ and $\Sigma(1:k,0:n)$ with cardinality
$\binom{n+1}{k+1}$ or $\binom{n+1}{k}$, respectively.
Furthermore, if $\sigma \in \Sigma(0:k,0:n)$, we denote by 
$\sigma^* \in \Sigma(1:n-k,0:n)$ the complementary 
map characterized by 
\begin{equation}\label{full-support}
\llbracket \sigma \rrbracket \cup \llbracket \sigma^* \rrbracket =
\{0,1,\ldots ,n\}.
\end{equation}
On the other hand, if $\sigma \in \Sigma(1:k,0:n)$, then $\sigma^*
\in \Sigma(0:n-k,0:n)$ is the complementary map such that
\eqref{full-support}
holds.

We will use the multi-index
notation $\alpha\in\mathbb N_0^n$, meaning $\alpha=(\alpha_1, \cdots,
\alpha_n)$ with integer $\alpha_i\ge 0$. We define $x^\alpha=x_1^{\alpha_1}\cdots x_n^{\alpha_n}$, and
$|\alpha|:=\sum\alpha_i$. We will also use the set $\mathbb N_0^{0:n}$ of
multi-indices $\alpha=(\alpha_0, \cdots, \alpha_n)$, with
$x^\alpha:=x_0^{\alpha_0}\cdots x_n^{\alpha_n}$.  The support $\llbracket\alpha\rrbracket$ of a multi-index $\alpha$ is
$\{\,i\,|\,\alpha_i>0\,\}$.
It is also useful to let
\[
\llbracket\alpha,\sigma\rrbracket =\llbracket\alpha\rrbracket 
\cup\llbracket\sigma\rrbracket , 
\quad \alpha\in\N_0^{0:n},\sigma\in\Sigma(j:k,l:m).
\]
If $\Omega \subset \R^n$ and $r \ge 0$, then $\P_r(\Omega)$ denotes the set of
real valued polynomials defined on $\Omega$ of degree less than or equal to
$r$.  For simplicity, we let $\P_r = \P_r(\R^n)$. Hence, if $\Omega$ has
nonempty interior, then $\dim \P_r(\Omega) = \dim \P_r = \binom{r+n}{n}$.
The case where $\Omega$ consists of a single point is allowed: then
$\P_r(\Omega)=\R$ for all $r\ge 0$.  For any $\Omega$, when
$r <0$, we take $\P_r(\Omega)=\{0\}$.

\subsection{Simplices and barycentric coordinates}
Let $T \in \R^n$ be an n-simplex with vertices $x_0,x_1,\ldots ,x_n$
in general position.  We let
$\Delta(T)$ denote all the subsimplices, or faces, of $T$, while $\Delta_k(T)$
denotes the set of subsimplices of dimension $k$. Hence, the cardinality of
$\Delta_k(T)$ is $\binom{n+1}{k+1}$.  We will use elements of the set
$\Sigma(j:k,0:n)$ to index the subsimplices of $T$.  For each
$\sigma\in\Sigma(j:k,0:n)$, we let $f_\sigma \in \Delta(T)$ be the closed
convex hull of the vertices $x_{\sigma(j)},\ldots ,x_{\sigma(k)}$, which we
henceforth denote by $[x_{\sigma(j)},\ldots ,x_{\sigma(k)}]$. Note that there
is a one-to-one correspondence between $\Delta_k(T)$ and $\Sigma(0:k,0:n)$. In
fact, the face $f_\sigma$ is uniquely determined by the range of $\sigma$,
$\llbracket \sigma \rrbracket$.  If $f = f_\sigma$ for $\sigma \in
\Sigma(j:k,0:n)$, we let the index set associated to $f$ be denoted by
$\I(f)$, i.e., $\I(f) = \llbracket \sigma \rrbracket$. If $f \in \Delta_k(T)$,
then $f^* \in \Delta_{n-k-1}(T)$ will denote the subsimplex of $T$ opposite
$f$, i.e., the subsimplex whose index set is the complement of
$\I(f)$ in $\{\, 0,1,\ldots ,n \, \}$. Note that if $\sigma \in \Sigma(0:k,0:n)$ and
$f=f_\sigma$, then $f^*= f_{\sigma^*}$.

We denote by $\lambda^T_0,\lambda^T_1,\dots,\lambda^T_n$ the barycentric
coordinate functions with respect to $T$, so $\lambda^T_i\in\P_1(T)$ is
determined by the equations $\lambda^T_i(x_j)=\delta_{ij}$, $0\le i,j\le n$.
The functions $\lambda^T_i$ form a basis for $\P_1(T)$, are non-negative on
$T$, and sum to $1$ identically on $T$.  Moreover, the subsimplices of $T$
correspond to the zero sets of the barycentric coordinates, i.e., if $f=
f_\sigma$ for $\sigma \in \Sigma(0:k,0:n)$, then $f$ is characterized by
\begin{equation*}
f = \{\,x\in T\,|\, \lambda^T_i(x)=0,\ i\in
\llbracket\sigma^*\rrbracket\,\}.
\end{equation*}
For a subsimplex $f \in \Delta(T)$, the barycentric coordinates
functions with respect to $f$, $\{\lambda_i^f\}_{i \in \I(f)}\subset \P_1(f)$, satisfy 
\begin{equation}\label{bary-restriction}
\lambda_i^f = \tr_{T,f} \lambda^T_i, \ i \in \I(f). 
\end{equation}
Here the
trace map $\tr_{T,f} : \P_1(T) \to \P_1(f)$ is the restriction of the
function to $f$. Due to the relation \eqref{bary-restriction},
we will sometimes omit the superscript $T$ or $f$, and simply write 
$\lambda_i$ instead of $\lambda_i^T$ or $\lambda_i^f$. Note
that, by linearity, the map $\lambda^f_i \to \lambda^T_i, \  i \in \I(f)$,
defines a \emph{barycentric extension operator}
$E^1_{f,T} : \P_1(f) \to \P_1(T)$, which is a right inverse of $\tr_{T,f}$.
The barycentric extension $E^1_{f,T}p$ can be characterized as the unique extension of the linear polynomial $p$ on $f$
to a linear polynomial on $T$ which vanishes on $f^*$.

\subsection{The Bernstein decomposition}\label{Bernstein}
Let $T = [x_0,x_1,\ldots ,x_n] \subset \R^n$ be as above and $ \{ \lambda_i
\}_{i=0}^n \subset \P_1(T)$ the corresponding barycentric coordinates. For $r
\ge 1$, the Bernstein basis for the space $\P_r(T)$ consists of
all monomials of degree $r$ in the variables $\lambda_i$, i.e., the basis
functions are given by
\begin{equation}\label{span1}
\{\,\lambda^\alpha = \lambda_0^{\alpha_0}
\lambda_1^{\alpha_1}\cdots \lambda_n^{\alpha_n}  \, 
| \, \alpha\in\mathbb N_0^{0:n}, \   |\alpha|=r \,\}.
\end{equation}
(It is common to take the scaled barycentric monomials
$(n!/\alpha!)\lambda^\alpha$ as the Bernstein basis elements, as in \cite{Lyche-Scherer},
but the scaling is not relevant here, and so we use the unscaled monomials.)
Of course, for $f \in \Delta(T)$, the space $\P_r(f)$ has the
corresponding basis
\[
\{\,(\lambda^f)^\alpha \, | \, \alpha\in\mathbb N_0^{0:n}, \ 
 |\alpha|=r, \  \llbracket \alpha \rrbracket \subseteq \I(f)\, \}.
\]
Hence, from this Bernstein basis, we also obtain a barycentric extension
operator,
$E=E^r_{f,T} : \P_r(f) \to \P_r(T)$,
by simply replacing $\lambda_i^f$ by $\lambda_i^T$ in the bases and
using linearity.

We let $\0\P_r(T)$ denote the subspace of
$\P_r(T)$ consisting of polynomials which vanish on the boundary of $T$ or,
equivalently, which are divisible by the corresponding bubble function
$\lambda_0\cdots\lambda_n$ on $T$.  Alternatively, we have 
\begin{equation}\label{span2}
\0\P_r(T)= \spn \{\lambda^\alpha \, | \,  \alpha\in\mathbb N_0^{0:n}, \ 
|\alpha|=r, \  \llbracket\alpha\rrbracket = \{0,\ldots ,n\}\,\}.
\end{equation}
Note that multiplication by the bubble function establishes an isomorphism
$\P_{r-n-1}(T)\cong\0\P_r(T)$.

\begin{figure}[htb]
\centerline{%
  \includegraphics[width=4in]{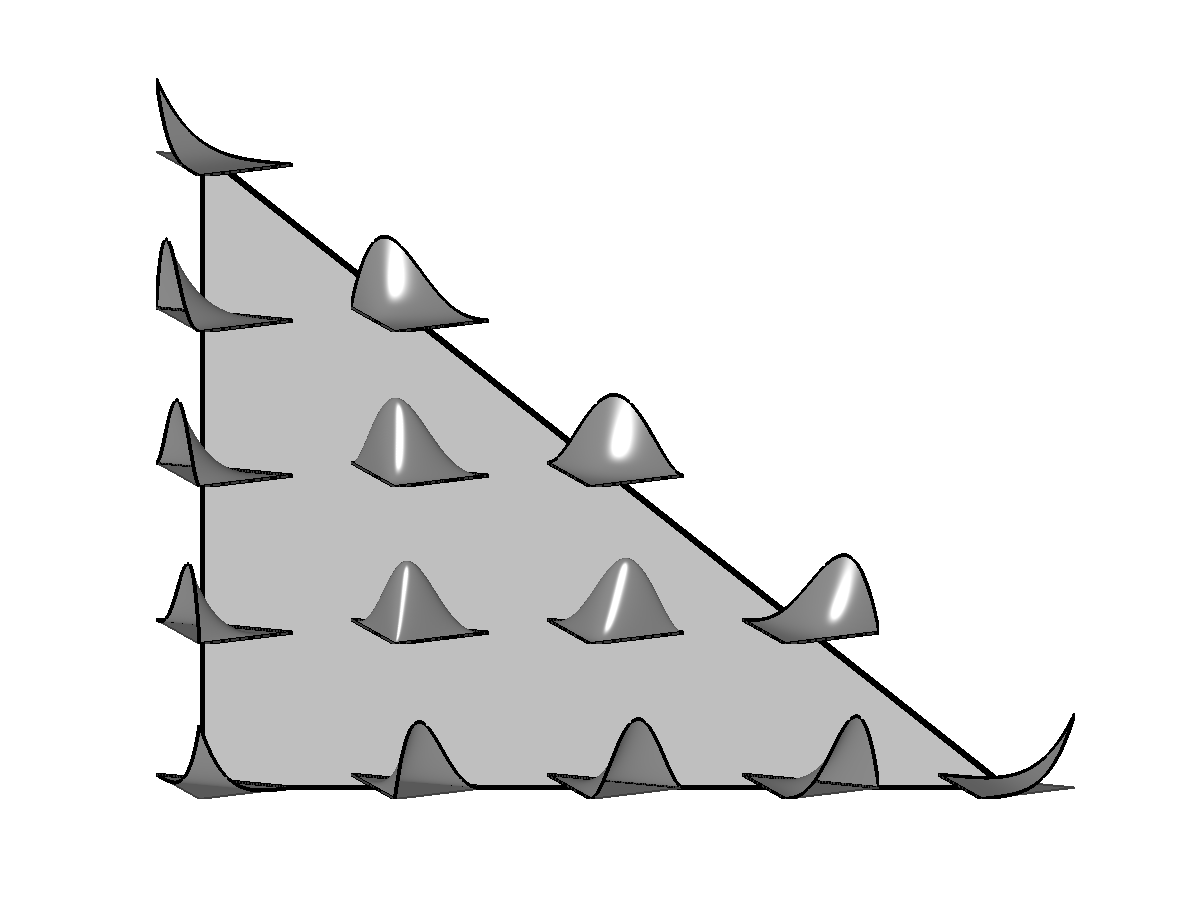}}
\vspace{-.4in}
\centerline{\quad$\underbrace{\hspace{1.75in}}$}
\centerline{\quad$\Big\uparrow$\rlap{$E$}}
\centerline{\quad$\overbrace{\hspace{1.75in}}$}
\centerline{\quad
  \includegraphics[width=1.75in]{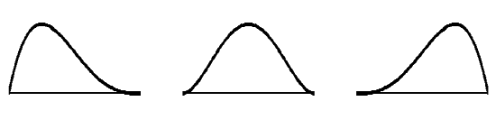}}
\vspace{-.2in}
\centerline{\rule{2.75in}{1pt}}
\caption[]{The Bernstein basis of $\P_4(T)$ for a triangle $T$.
One basis function is associated with each vertex, three with each edge,
and three with the triangle.  The basis functions associated with any
face $f$ are obtained by extending basis functions for $\0\P_4(f)$ to the triangle.}
\label{fg:Bernsteinbasis}
\end{figure}

The Bernstein basis \eqref{span1} leads to an explicit geometric decomposition
of the space $\P_r(T)$. Namely, we associate to the face $f$, the subspace of
$\P_r(T)$ that is  spanned by the  basis functions $\lambda^\alpha$ with
$\llbracket\alpha\rrbracket = \I(f)$.
We then note that this subspace is precisely $E[\0\P_r(f)]$, i.e.,
\begin{equation}\label{def-PTf0}
E[\0\P_f(f)] = \spn\{\,\lambda^\alpha\,|\,\alpha\in\mathbb N_0^{0:n},\  
|\alpha|=r,\  \llbracket \alpha \rrbracket = \I(f)\,\}.
\end{equation}
Clearly,
\begin{equation}\label{B-decomp}
\P_r(T) = \bigoplus_{f\in\Delta(T)} E[\0\P_r(f)],
\end{equation}
which we refer to as the \emph{Bernstein
decomposition} of the space $\P_r(T)$.
This is an example of a geometric decomposition, as discussed
in the introduction.
%Moreover, the decomposition \eqref{B-decomp} is consistent
%with the decomposition \eqref{dual-decomp} of the dual space in the
%sense that the restriction of the elements of $W_r(T,f)$ to $\0\P_r(T,f)$
%gives an isomorphism $W_r(T,f)\cong \0\P_r(T,f)^*$.
%Note also that for any multi-index $\alpha \in N_0^{0:n}$,
%the corresponding polynomial $\lambda^\alpha$ satisfies
%\[
%\tr_{T,f} \lambda^\alpha = 0 \quad \text{if and only if } \llbracket
%\alpha \rrbracket \cap \I(f^*) \neq \emptyset.
%\]
%As a consequence, if $\omega \in \P_r(T,f)$ and $g \in \Delta(T)$, then the
%trace of $\omega$ 
%on $g$, $\tr_{T,g}\omega$, satisfies
%\[
%\tr_{T,g}\omega = 0, \quad \text{if } \dim g \le \dim f, \, g \neq f.
%\]
%From this property and a simple inductive argument, it can easily be seen that
%the sum in \eqref{B-decomp} is direct.
An illustration of the decomposition
\eqref{B-decomp} is given in Figure~\ref{fg:Bernsteinbasis}.

Moreover, the extension operator $E$ may also be characterized
geometrically, without recourse to barycentric coordinates.
To obtain such a characterization,
we first recall that a smooth function $u : T \to \R $ is said to {\em
vanish to order $r$} at a point $x$ if
\begin{equation*}
(\partial^\alpha u)(x) = 0, \  \alpha \in\mathbb N_0^n, \  |\alpha| \le r-1.
\end{equation*}
We also say that $u$ vanishes to order $r$ on a set $g$ if it vanishes to order
$r$ at each point of $g$.
Note that the extension operator $E=E_{f,T}^r$ has the property that
for any $\mu \in \P_r(f)$, $E \mu$ vanishes to order $r$ on $f^*$.  In
fact, if we set
\begin{equation*}
\P_r(T,f) = \{\,\omega \in \P_r(T) \,|\,
\omega \text{ \rm vanishes to order } r \text{ \rm on } f^*\, \},
\end{equation*}
we can prove

\begin{lem}\label{0-forms}
$\P_r(T,f) = E[\P_r(f)]$ and for $\mu \in \P_r(f)$, $E \mu =E_{f,T}^r
\mu$ can be characterized as the unique extension of $\mu$ to
$\P_r(T,f)$.
\end{lem}
\begin{proof}
It is easy to see that $E[\P_r(f)] \subseteq \P_r(T,f)$. To establish
the reverse inclusion, we observe that if $f^* = \{x_i\}$, then $\omega
\in \P_r(T)$
vanishes at $f^*$ if and only if it can be written in the form
\begin{equation*}
\omega = \sum_{|\alpha|=r} c_{\alpha} \lambda^{\alpha}, 
\end{equation*}
where the sum is restricted to multi-indices $i$ for
which $\alpha_i=0$.  For a more general set 
$f^*$, this fact will be true for any $i \in \I(f^*)$ and hence
\begin{equation*}
\omega = \sum_{\substack{|\alpha|=r \\ \llbracket\alpha\rrbracket 
\subseteq \I(f)}}
 c_{\alpha} \lambda^{\alpha}, 
\end{equation*}
and so $\omega \in E[\P_r(f)]$.
\end{proof}

We also note that it follows immediately from Lemma~\ref{0-forms}
that the space $E[\0 \P_r(f)]$ appearing in the Bernstein decomposition
\eqref{B-decomp} is characterized by
\begin{equation*}
E[\0 \P_r(f)] = \{\omega \in \P_r(T) \,|\,
\omega \text{ vanishes to order } r \text{ on } f^*, \, 
\tr_{T,f} \omega \in \0 \P_r(f)\}.
\end{equation*}

In this paper, we will establish results analogous to those of this section
for spaces of polynomial differential forms, and in particular, direct sum
decompositions of these spaces analogous to the Bernstein decomposition
\eqref{B-decomp}.

\subsection{Differential forms}\label{df}
Next we indicate the notations we will be using for basic concepts related to
differential forms.  See \cite[\S 2]{acta} or the references indicated
there for a more detailed treatment.  For $k\ge0$, we denote by $\Alt^k V$ the
set of real-valued, alternating $k$-linear maps on a vector space $V$ (with
$\Alt^0V=\R$). 
Hence, $\Alt^kV$ is a vector space of dimension $\binom{\dim
V}{k}$.  The exterior product, or the wedge product, maps $\Alt^j V\times
\Alt^k V$ into $\Alt^{j+k} V$.  If $\omega \in \Alt^k V$ and $v \in V$, then
the contraction of $\omega$ with $v$, $\omega \lrcorner v \in \Alt^{k-1} V$,
is given by $\omega \lrcorner v(v_1, \ldots,v_{k-1}) = \omega(v,v_1,
\ldots,v_{k-1})$.

If $\Omega$ is a smooth manifold (e.g., an open subset of Euclidean
space), a differential $k$-form on $\Omega$ is a map which
assigns to each $x\in\Omega$ an element of $\Alt^kT_x\Omega$, where
$T_x\Omega$ is the tangent space to $\Omega$ at $x$.
In case $f$ is an open subset of an affine subspace of Euclidean, all the
tangents spaces $T_x f$ may be canonically identified, and we simply write
them as $T_f$. 

We denote by $\Lambda^k(\Omega)$ the
space of all smooth differential $k$-forms on $\Omega$.  The exterior
derivative $d$ maps $\Lambda^k(\Omega)$ to $\Lambda^{k+1}(\Omega)$.
It satisfies $d\circ d=0$, so defines a complex
\[
0\to\Lambda^0(\Omega)
\xrightarrow{d} \Lambda^1(\Omega) \xrightarrow{d} 
\cdots \xrightarrow{d}
\Lambda^n(\Omega)\rightarrow{}0,
\]
the de~Rham complex.
%We recall the Leibniz rule for the exterior derivative:
%\begin{equation}\label{leibniz}
%d(\omega \wedge \eta) = d \omega \wedge \eta + (-1)^j\omega \wedge d
%\eta,
%\quad \omega \in \Lambda^j(\Omega), \eta \in \Lambda^k(\Omega).
%\end{equation}

If $F:\Omega\to\Omega'$, is a smooth map between smooth manifolds, 
then the pullback 
$F^* :\Lambda^k(\Omega') \to \Lambda^k(\Omega)$ is given by 
\[
(F^*\omega)_x(v_1, v_2,\ldots ,v_k) =
\omega_{F(x)}(DF_x(v_1),
DF_x(v_2), \ldots ,DF_x(v_k)),
\]
where the linear map $D F_x : T_x\Omega \to T_{F(x)}\Omega'$ is the derivative
of $F$ at $x$.
The pullback commutes with the exterior derivative, i.e.,
\[
F^*(d \omega) = d(F^*\omega), \quad \omega \in
\Lambda^k(\Omega'),
\]
and distributes with respect to the wedge product:
\[
F^*(\omega \wedge \eta) = F^* \omega \wedge F^* \eta.
\]
We also recall the integral of a $k$-form over an orientable $k$-dimensional
manifold is defined, and
\begin{equation}\label{int-pull-back}
\int_\Omega F^*\omega = \int_{\Omega'} \omega,
\quad \omega \in \Lambda^n(\Omega'),
\end{equation}
when $F$ is an orientation-preserving diffeomorphism.

If $\Omega'$ is a submanifold of $\Omega$, then the pullback of the
inclusion $\Omega'\hookrightarrow\Omega$ is the \emph{trace map}
$\tr_{\Omega,\Omega'}:\Lambda^k(\Omega)\to\Lambda^k(\Omega')$. 
If the domain $\Omega$ is clear from the context, we may write
$\tr_{\Omega'}$ instead of $\tr_{\Omega,\Omega'}$,
and if $\Omega'$ is the boundary of $\Omega$, $\partial \Omega$, we
just write $\tr$.
Note that if $\Omega'$ is a submanifold
of positive codimension and $k>0$, then the vanishing
of $\Tr_{\Omega,\Omega'}\omega$ on $\Omega'$ for $\omega\in\Lambda^k(\Omega)$
does \emph{not} imply that $\omega_x\in\Alt^kT_x\Omega$ vanishes
for $x\in\Omega'$, only that it vanishes when applied to $k$-tuples
of vectors tangent to $\Omega'$, or, in other words, that the tangential
part of $\omega_x$ with respect to $T_x\Omega'$ vanishes.

If $\Omega$ is a subset of $\R^n$ (or, more generally, a Riemannian
manifold), we can define 
the Hilbert space $L^2\Lambda^k(\Omega) \supset \Lambda^k(\Omega)$
of $L^2$ differential $k$-forms,
and the Sobolev space
\[
H\Lambda^k(\Omega) := \{\, \omega \in L^2\Lambda^k(\Omega) \, | \, 
d\omega \in L^2\Lambda^{k+1}(\Omega)\, \}.
\]
The $L^2$ de Rham complex  is the sequence of mappings and spaces given by
\begin{equation}\label{deRham}
0\to H\Lambda^0(\Omega)
\xrightarrow{d} H\Lambda^1(\Omega) \xrightarrow{d} 
\cdots \xrightarrow{d}
H\Lambda^n(\Omega)\to0.
\end{equation}
We remark that for $\Omega\subset\R^n$, $H\Lambda^0(\Omega)$ is
equal to the ordinary Sobolev space
$H^1(\Omega)$ and, via the identification of $\Alt^n\R^n$ with $\R$,
$H\Lambda^n(\Omega)$ can be identified with
$L^2(\Omega)$.  Furthermore, in the case $n=3$, the spaces $\Alt^1 \R^3$
and $\Alt^2 \R^3$ can be identified with $\R^3$, and the complex
\eqref{deRham} may be identified with the complex
\[
0\to H^1(\Omega)
\xrightarrow{\grad} H(\curl;\Omega) \xrightarrow{\curl} 
H(\div;\Omega) \xrightarrow{\div} 
L^2(\Omega)\to0.
\]

\subsection{Simplicial triangulations}\label{st}
Let $\Omega$ be a bounded polyhedral domain in $\R^n$ and $\T$
a finite set of $n$-simplices. We will refer to $\T$
as a \emph{simplicial triangulation} of $\Omega$ if the union of all 
the elements of $\T$ is the closure of $\Omega$, and the intersection
of two is either empty or a common subsimplex of each. For $0 \le j
\le n$, we let 
\[
\Delta_j(\T) = \bigcup_{T \in \T}\Delta_j(T) \quad \text{and }
\Delta(\T)= \bigcup_{j=0}^n \Delta_j(\T).
\]
In the finite element exterior calculus,  
we employ spaces of differential forms $\omega$
which are piecewise smooth (usually polynomials) with respect to $\T$,
i.e., the restriction $\omega|_T$ is smooth for each $T \in \T$. Then for $f\in\Delta_j(\T)$ with $j\ge k$,
$\tr_f\omega$ may be multi-valued, in that we can assign
a value for each $T\in\T$ containing $f$ by first restricting
$\omega$ to $T$ and then taking the trace on $f$.  If all such
traces coincide, we say that $\tr_f\omega$ is single-valued.
The following
lemma, a simple consequence of Stokes' theorem,
cf.~\cite[Lemma 5.1]{acta}, is a key result.

\begin{lem}\label{continuity}
Let $\omega \in L^2\Lambda^k(\Omega)$ be piecewise smooth 
with respect to the triangulation $\mathcal T$. The following
statements are equivalent:
\begin{itemize}
\item[(1)] $\omega \in H\Lambda^k(\Omega)$,
\item[(2)] $\tr_f \omega$ is single-valued for all $f \in
  \Delta_{n-1}(\T)$,
\item[(3)] $\tr_f \omega$ is single-valued for all $f \in \Delta_j(\T)$, $k
  \le j \le n-1$.
\end{itemize}
\end{lem}
As a  consequence of this theorem,
in order to construct subspaces of $H\Lambda^k(\Omega)$, consisting of
differential forms $\omega$ which are piecewise smooth with respect to
the 
triangulation $\T$, we need to build into the construction that $\tr_f \omega$
is single-valued for each $f \in \Delta_j(\T)$ for $k \le j \le n-1$.

\section{Polynomial and piecewise polynomial differential forms}\label{sec:pdf}
In this section we formally define the two families of spaces of polynomial differential
forms $\P_r\Lambda^k$ and $\P^-_r\Lambda^k$.
These polynomial spaces will then be used to define piecewise 
polynomial differential forms with respect to a simplicial triangulation
of a bounded polyhedral domain in $\R^n$. In fact, as explained in
\cite[\S 3.4]{acta}, the two families presented here are nearly the
only affine invariant spaces of polynomial differential forms.

\subsection{The space $\P_r \Lambda^k$}
Let $\Omega$ be a subset of $\R^n$.
For $0\le k\le n$, we
let $\P_r\Lambda^k(\Omega)$ be the subspace of $\Lambda^k(\Omega)$
consisting of all $\omega \in \Lambda^k(\Omega)$ such that 
$\omega(v_1,v_2,\ldots ,v_k) \in \P_r(\Omega)$ for each choice of
vectors $v_1,v_2,\ldots v_k \in \R^n$. Frequently, we will write 
$\P_r\Lambda^k$ instead of $\P_r\Lambda^k(\R^n)$.
The space $\P_r\Lambda^k$ is isomorphic to $\P_r\otimes \Alt^k$ and
\begin{equation}\label{pdim}
\dim\P_r\Lambda^k = \dim\P_r \x \dim\Alt^k\R^n = \binom{r+n}{n}\binom{n}{k}
= \binom{r+k}{r}\binom{n+r}{n-k}.
\end{equation}
Furthermore, if $\Omega \subset \R^n$ with nonempty interior, then 
$\dim\P_r\Lambda^k(\Omega) \cong \dim\P_r\Lambda^k$.

If $T$ is a simplex, we define
\begin{equation*}
\0\P_r\Lambda^k(T) = \{\,\omega\in\P_r\Lambda^k(T)\,|\,
\tr\omega=0\,\}.
\end{equation*}
In the case $k=0$, this space simply consists of all the
polynomials divisible by the bubble function
$\lambda_0\cdots\lambda_n$, so
\begin{equation}\label{a}
\0\P_r\Lambda^0(T) \cong\P_{r-n-1}\Lambda^0(T).
\end{equation}
For $k=n$, the trace map
vanishes, so we have
\begin{equation}\label{b}
\0\P_r\Lambda^n(T)=\P_r\Lambda^n(T).
\end{equation}

\subsection{The space $\P^-_r\Lambda^k$}
The Koszul differential $\kappa$ of a differential $k$-form $\omega$
on $\R^n$
is the $(k-1)$-form given by
\begin{equation*}
(\kappa\omega)_x(v_1,\ldots,v_{k-1}) =
\omega_x\bigl(X(x),v_1,\ldots,v_{k-1}\bigr),
\end{equation*}
where $X(x)$ is the vector from the origin to $x$.
For each $r$, $\kappa$
maps $\P_{r-1}\Lambda^{k}$ to $\P_r\Lambda^{k-1}$, and the
Koszul complex
\[
\begin{CD}
0\to \P_{r-n}\Lambda^n  @>\kappa>>
\P_{r-n+1}\Lambda^{n-1} @>\kappa>> \cdots  @>\kappa>> \P_r\Lambda^0
\to \R \to0,
\end{CD}
\]
is exact. Furthermore, the Koszul operator satisfies the Leibniz
relation
\begin{equation}\label{kappa-leibniz}
\kappa(\omega \wedge \eta) 
= (\kappa\omega)\wedge \eta + (-1)^k\omega \wedge (\kappa \eta), \quad
\omega \in \Lambda^k, \  \eta \in \Lambda^l.
\end{equation}
We define
\begin{equation*}
\P^-_r\Lambda^k = \P^-_r\Lambda^k(\R^n) =\P_{r-1}\Lambda^k+\kappa\P_{r-1}\Lambda^{k+1}.
\end{equation*}
From this definition, we easily see that $\P^-_r\Lambda^0=\P_r\Lambda^0$ and
$\P^-_r\Lambda^n=\P_{r-1}\Lambda^n$. However, if $0<k<n$, then 
\begin{equation*}
\P_{r-1}\Lambda^k\subsetneq \P^-_r\Lambda^k
\subsetneq \P_r\Lambda^k.
\end{equation*}
An important property of the spaces $\P^-_r\Lambda^k$
is the closure relation
\begin{equation}\label{closure}
\P^-_r\Lambda^k \wedge \P_s^-\Lambda^l \subseteq
\P_{r+s}^-\Lambda^{k+l}.
\end{equation}
A key identity relating the Koszul operator $\kappa$ with the 
exterior derivative $d$ is the homotopy relation
\begin{equation}\label{dk+kd}
(\D\kappa+\kappa \D)\omega = (r+k)\omega,\quad
\omega\in\H_r\Lambda^k,
\end{equation}
where $\H_r\Lambda^k$ is the space of \emph{homogeneous} polynomial
$k$-forms of degree $r$.

Using the homotopy relation and the exactness of the
Koszul complex, we can inductively
compute the dimension of $\P^-_r\Lambda^k$ as
\begin{equation}\label{pminusdim}
\dim \P^-_r\Lambda^k = \binom{r+k-1}{k}\binom{n+r}{n-k}.
\end{equation}
If $\Omega \subset \R^n$, then $\P^-_r\Lambda^k(\Omega)$ denotes the
restriction of functions in $\P^-_r\Lambda^k$ to $\Omega$, which
implies that the space $\P^-_r\Lambda^k(\Omega)$ is isomorphic to 
$\P^-_r\Lambda^k$ if $\Omega$ has nonempty interior.
Finally, we remark that although the Koszul operator $\kappa$ depends on the
choice of origin used to associate a point in $\R^n$ with a vector, the space
$\P^-_r\Lambda^k$ is unaffected by the choice of origin.
We refer to \cite{acta} for more details on the spaces $\P^-_r\Lambda^k$.
In particular, if $T \subset \R^n$ is a simplex and $f \in \Delta_j(T)$ then $\tr_f
\P_r^-\Lambda^k(T)
= \P_r^-\Lambda^k(f)$, where the space
$\P_r^-\Lambda^k(f)\cong\P_r^-\Lambda^k(\R^j)$ depends on $f$,
but is independent of $T$. 

For a simplex $T$, we define
\begin{equation*}
\0\P^-_r\Lambda^k(T) = \{\,\omega\in\P^-_r\Lambda^k(T)\,|\,
\tr\omega=0\,\}.
\end{equation*}

From the Hodge star isomorphism, we have that $\P_{r-n-1}\Lambda^0(T)\cong\P_{r-n-1}\Lambda^n(T)=\P_{r-n}^-\Lambda^n(T)$
and that $\P_r\Lambda^n(T)\cong\P_r\Lambda^0(T)=\P^-_r\Lambda^0(T)$.  Therefore
\eqref{a} and \eqref{b} become
\begin{equation}\label{c}
\0\P_r\Lambda^0(T) \cong\P_{r-n}^-\Lambda^n(T), \quad \0\P_r\Lambda^n(T)\cong\P^-_r\Lambda^0(T).
\end{equation}
These are the two extreme cases of the relation
\begin{equation}\label{iso0}
\0\P_r\Lambda^k(T) \cong \P_{r-n+k}^-\Lambda^{n-k}(T), \quad 0\le k\le n.
\end{equation}
But \eqref{c} can also be written
\begin{equation*}
\0\P^-_r\Lambda^0(T) \cong\P_{r-n-1}\Lambda^n(T), \quad \0\P^-_r\Lambda^n(T)\cong\P_{r-1}\Lambda^0(T),
\end{equation*}
(where we have substituted $r-1$ for $r$ in the second relation), which are the extreme cases of
\begin{equation}\label{iso1}
\0\P^-_r\Lambda^k(T) \cong \P_{r-n+k-1}\Lambda^{n-k}(T), \quad 0\le k\le n.
\end{equation}
That the isomorphisms in \eqref{iso0} and \eqref{iso1} do indeed exist for all
$k$ follows from Corollary~\ref{iso-spaces} below.

\subsection{The spaces $\P_r\Lambda^k(\T)$ and $\P^-_r\Lambda^k(\T)$ }
For $\T$ a simplicial triangulation of a domain $\Omega \in
\R^n$, we define
\begin{multline*}
\P_r\Lambda^k(\T) =\{\,\omega\in L^2\Lambda^k(\Omega)\,|\, 
\omega|_T\in\P_r\Lambda^k(T) \quad \forall T\in\T,
\\*[-3pt] 
\tr_f \omega  \text{ is single-valued for }f \in
\Delta_j(\T),\  k \le j \le n-1 \,\},
\end{multline*}
and define $\P^-_r\Lambda^k(\T)$ similarly.
In view of Lemma~\ref{continuity}, we have
\begin{align*}
\P_r\Lambda^k(\T)&=\{\,\omega\in H\Lambda^k(\Omega)\,|\, 
\omega|_T\in\P_r\Lambda^k(T) \quad \forall T\in\T\,\},
\\
\P^-_r\Lambda^k(\T)&=\{\,\omega\in H\Lambda^k(\Omega)\,|\, 
\omega|_T\in\P^-_r\Lambda^k(T) \quad \forall T\in\T\,\}.
\end{align*}

\section{Consistent extension operators and geometric decompositions}\label{cgd}

Let $\T$ be a simplicial triangulation of $\Omega \subset \R^n$, and
let there be given a finite
dimensional subspace $X(T)$ of $\Lambda^k(T)$
for each $T\in\T$.  In this section we shall define the notion of
a \emph{consistent family of extension operators}, and show that it
leads to the construction of a geometric decomposition  and a local basis of the
finite element space
\begin{equation}\label{defxt}
X(\T) = \{\,\omega\in L^2\Lambda^k(\Omega)\,|\,
\omega|_T\in X(T) \ \forall T\in\T,\
\tr_f \omega  \text{ is single-valued for }f \in
\Delta(\T)\,\}. 
\end{equation}
We note that as a result of Lemma~\ref{continuity}, $X(\T) \subset
H\Lambda^k(\Omega)$.

For the Lagrange finite element space $\P_r(\T)=\P_r\Lambda^0(\T)$, both the Bernstein
basis discussed in Section~\ref{Bernstein} and the dual basis
discussed in the next section arise from this construction. One of the main goals of this paper is to generalize these bases to the two families of finite element
spaces of $k$-forms.

We require that the family of spaces $X(T)$ fulfills the following consistency assumption:
\begin{equation}\label{cassump}
\tr_{T,f}X(T)=\tr_{T',f}X(T') \text{ whenever $T,T'\in\T$ with $f\in\Delta(T)\cap\Delta(T')$}.
\end{equation}
In this case, we may define for any $f\in\Delta(\T)$,
$X(f)=\tr_{T,f}X(T)$ where $T\in\T$ is any simplex containing $f$.
We also define $\0X(f)$ as the subspace of $X(f)$
consisting of all $\omega \in X(f)$ such that $\tr_{f,\partial f} \omega = 0$.
Note that
\begin{equation}\label{traceprop}
\tr_{g,f}X(g) =X(f) \text{ for all $f,g\in\Delta(\T)$ with $f\subseteq g$}.
\end{equation}
Consequently, for each such $f$ and $g$ we may choose an
\emph{extension operator}
$E_{f,g}:X(f)\to X(g)$, i.e., a right inverse of 
$\tr_{g,f}:X(g) \to X(f)$.
%Let $X(g,f)\subset X(g)$ denote
%the range of the extension operator.
We say that a family of extension operators $E_{f,g}$, defined for
all $f,g\in\Delta(\T)$ with $f\subseteq g$, is \emph{consistent} if
%\begin{equation}\label{H1}
%X(g,e)\subset X(g,f)  \text{ for all $e,f,g\in\Delta(T)$ with $e\subset %f\subset g$},
%\end{equation}
\begin{equation}\label{H2}
\tr_{h,g} E_{f,h}=E_{f\cap g,g}\tr_{f,f\cap g} \text{ for all $f,g,h\in\Delta(\T)$ with $f,g\subseteq h$}.
\end{equation}
In other words, we require that the diagram
\begin{equation*}
\begin{CD}
X(f) @>E>> X(h) \\
@VV\tr V @VV \tr V \\
X(f\cap g) @>E>> X(g)
\end{CD}
\end{equation*}
commutes. 

One immediate implication of \eqref{H2} is that
for $\omega \in X(f)$,
\begin{equation}\label{consistent}
\tr_{h,g} E_{f,h} \omega = E_{f,g} \omega
\text{\quad for all $f,g,h\in\Delta(\T)$ with
$f \subseteq g \subseteq h$}.
\end{equation}
A second implication is:
\begin{lem}\label{zero-trace-prop}
Let $h\in\Delta(\T)$, and $f,g \in \Delta(h)$ with $f \nsubseteq g$.
Then $\tr_{h,g}\omega = 0$ for all $\omega \in E_{f,h}\0X(f)$.
\end{lem}
\begin{proof} Let $\omega =E_{f,h}\mu$ with $\mu\in \0X(f)$. 
Since $f \nsubseteq g$, we have 
$f \cap g \subset \partial f$, and therefore
$\tr_{f,f \cap g}\mu = 0$. Then, by \eqref{H2},
$\tr_{h,g}\omega= \tr_{h,g}E_{f,h}\mu
=E_{f\cap g,g}\tr_{f,f\cap g}\mu
=0$.
\end{proof}

We now define an extension operator $E_f:\0X(f)\to X(\T)$ for each $f\in\Delta(\T)$.
Given $\mu\in\0X(f)$, we define $E_f\mu$ piecewise:
\begin{equation}\label{eft}
(E_f \mu)|_T=
\begin{cases}
 E_{f,T} \mu &\text{ if $f\subseteq T$},
\\
0, &\text{ otherwise}.
\end{cases}
\end{equation}
We claim that for each $g\in\Delta(\T)$, $\tr_g E_f\mu$ is single-valued,
so $E_f\mu$ does indeed belong to $X(\T)$.  To see this, we consider separately
the cases $f\subseteq g$ and $f\nsubseteq g$.  In the former case, if
$T\in\T$ is any simplex containing $g$, then $f\subseteq T$, and so
\begin{equation*}
 \tr_{T,g}[(E_f \mu)|_T]=\tr_{T,g}E_{f,T} \mu = E_{f,g}\mu
\end{equation*}
by \eqref{consistent}.  Thus $\tr_{T,g}[(E_f \mu)|_T]$ does not depend
on the choice of $T$ containing $g$, so in this case we have established
that $\tr_g E_f\mu$ is single-valued.  On the other hand, if $f\nsubseteq g$
then $\tr_{T,g}[(E_f \mu)|_T]=0$ for any $T$ containing $g$, either because
$f\nsubseteq T$ and so $(E_f \mu)|_T=0$,  or by Lemma~\ref{zero-trace-prop}
if $f\subseteq T$.  Thus we have established that all traces of $E_f\mu$ are single-valued, and so we have defined
extension operators $E_f:\0X(f)\to X(\T)$ for each $f\in\Delta(\T)$.
We refer to $E_f$ as the global extension operator determined by the consistent family of
extension operators.

We easily obtain this variant of Lemma~\ref{zero-trace-prop}.
\begin{lem}\label{zero-trace-mesh}
Let $f,g\in\Delta(\T)$, $f\nsubseteq g$.  Then
$\tr_g\omega=0$ for all $\omega\in E_f\0X(f)$.
\end{lem}
\begin{proof}
Pick $T\in\T$ containing $g$.  If $f\in\Delta(T)$, then we can apply
Lemma~\ref{zero-trace-prop} with $h=T$.  Otherwise, $(E_f\mu)|_T=0$
for all $\mu\in \0X(f)$.
\end{proof}
The following theorem is the main result of this section.
\begin{thm}\label{basic-decomp}
Let $\T$ be a simplicial triangulation and suppose that for each $T\in\T$, a
finite-dimensional subspace $X(T)$ of $\Lambda^k(T)$ is given fulfilling the consistency
assumption \eqref{cassump}.  Assume that there is a consistent family
of extensions operators $E_{g,f}$ for all $f,g\in\Delta(\T)$ with $f\subseteq g$.  Define
$E_f$, $f\in\Delta(\T)$ by \eqref{eft}.  Then the space
$X(\T)$ defined in \eqref{defxt} admits the direct sum decomposition
\begin{equation}\label{x-decomp-mesh}
 X(\T) =  \bigoplus_{f\in\Delta(\T)} E_f\0X(f).
\end{equation}
\end{thm}

\begin{proof}
To show that the sum is direct, we assume that $\sum_{f \in \Delta(\T)}
\omega_f = 0$, where $\omega_f \in E_{f,T}\0X(f)$, and prove by
induction that $\omega_f = 0$ for all $f \in \Delta(\T)$ with $\dim f
\le j$.  This is certainly true for $j<k$, (since then $\Lambda^k(f)$
and, {\it a fortiori}, $X(f)$ vanishes), so we assume it is true and must
show that $\omega_g=0$ for $g \in \Delta_{j+1}(\T)$.  By
Lemma~\ref{zero-trace-mesh},
\[
0 = \tr_g\Big(\sum_{f\in\Delta(\T)} \omega_f\Big) = \tr_g\omega_g.
\]
Hence, $\omega_g = E_g\tr_g\omega_g=0$. We thus conclude that the
sum is direct, and 
$X(\T) \supseteq \bigoplus_{f\in\Delta(\T)} E_f\0X(f)$.

To show that this is an equality, we write any $\omega \in X(\T)$ 
in the form
\[
\omega = \omega^n - \sum_{j=k}^{n-1}(\omega^{j+1} - \omega^j),
\]
where $\omega^k = \omega$, and  for $k < j \le n$, 
$\omega^j \in X(\T)$ is defined recursively by
\[
\omega^{j+1} = \omega^j - \sum_{f \in \Delta_j(\T)} E_f
\tr_f\omega^j.
\]
We shall prove by induction that for $k \le j \le n$
\begin{equation}\label{ind}
\tr_f\omega^j \in \0X(f), \quad f \in \Delta_j(\T).
\end{equation} 
Assuming this momentarily, we get that
$\omega^{j+1} - \omega^j \in \sum_{f \in\Delta_j(\T)} E_f\0X(f)$.
Also, $\omega^n|_T=\tr_T\omega^n\in\0X(T)$ for all $T\in\T$,
and  $\omega^n=\sum_{T\in\T}\tr_T(\omega^n|_T)$.  Thus,
$\omega \in \bigoplus_{f\in\Delta(T)} E_{f,T}\0X(f)$ as desired.

To prove \eqref{ind} inductively, we first note it is certainly true if
$j=k$, since $\0X(f)=X(f)$ for $f\in \Delta_k(\T)$.  Now
assume \eqref{ind}
and let $g \in \Delta_{j+1}(\T)$.   We show that
$\tr_g\omega^{j+1} \in \0X(g)$, by showing that $\tr_h\omega^{j+1} =0$
for $h \in \Delta_j(g)$.  In fact,
\[
\tr_h\omega^{j+1} = \tr_h \omega^j 
- \sum_{f \in \Delta_j(\T)} \tr_h E_f \tr_f\omega^j.
\]
Now $\tr_f\omega^j \in \0X(f)$ by the inductive hypothesis,
and therefore,  by Lemma~\ref{zero-trace-mesh},
$\tr_h E_f \tr_f\omega^j =0$ unless $f=h$, in which case
$\tr_h E_f \tr_f\omega^j =\tr_h\omega^j$. Thus,
\[
\tr_h\omega^{j+1} = \tr_h \omega^j 
- \tr_h \omega^j =0.
\]
This completes the proof of the theorem.
\end{proof}

\begin{remark}
By considering the case of a mesh consisting of a single simplex $T$, we see that
\begin{equation}\label{x-decomp}
 X(T) =  \bigoplus_{f\in\Delta(T)} E_{f,T}\0X(f).
\end{equation}
\end{remark}

The decomposition \eqref{x-decomp-mesh} is very important in
practice. It leads immediately to a local basis for the large space
$X(\T)$ consisting of elements $E_f\mu$, where $f$ ranges over $\Delta(\T)$  and $\mu$ ranges
over a basis for the space $\0X(f)$.

We close this section with the simplest example of this theory.
Let $X(T) = \P_r(T)=\P_r\Lambda^0(T)$ be the polynomial space discussed in
Section~\ref{Bernstein}.  Then \eqref{cassump} is fulfilled and the trace spaces $X(f)$ are simply
$\P_r(f)$ for $f\in\Delta(\T)$.  For $f,g\in\Delta(\T)$ with $f\subseteq g$,
the trace operator $\tr_{g,f}$ and barycentric extension operator
$E_{f,g}$ are given in barycentric coordinates as follows.  If
$\alpha\in\N^{0:n}_0$ with $\llbracket \alpha\rrbracket\subseteq\I(g)$, then
\begin{equation*}
 \tr_{g,f}(\lambda^g)^\alpha=
\begin{cases}
 (\lambda^f)^\alpha & \text{if $\llbracket g\rrbracket\subseteq\I(f)$},
\\
0, &\text{otherwise}.
\end{cases}
\end{equation*}
For $\alpha\in\N^{0:n}_0$ with $|\alpha|=r$ and $\llbracket \alpha\rrbracket\subseteq\I(f)$,  then
$E_{f,g}(\lambda^f)^\alpha=(\lambda^g)^\alpha$.  We now check that
the family of barycentric
extension operators is consistent, i.e., we verify \eqref{H2}.
We must show that if $f,g,h\in\Delta(\T)$ with $f,g\subseteq h$, then
\begin{equation*}
 \tr_{h,g}E_{f,h}(\lambda^f)^\alpha
= E_{f\cap g,g}\tr_{f,f\cap g}(\lambda^f)^\alpha
\end{equation*}
for all multi-indices $\alpha$ with $|\alpha|=r$ and
$\llbracket\alpha\rrbracket\subseteq\I(f)$.  Indeed, it is easy to
check that both sides are equal $(\lambda^g)^\alpha$ if
$\llbracket\alpha\rrbracket\subseteq\I(g)$ and zero otherwise.  Note
that, in this case, the decomposition \eqref{x-decomp} is simply the
Bernstein decomposition \eqref{B-decomp}.  If we then define, as in
the general definition \eqref{defxt} above,
\begin{equation*}
\P_r \Lambda^0(\T) = \{\,\omega\in L^2(\Omega)\,|\,
\omega|_T\in \P_r(T) \ \forall T\in\T,\
\tr_f \omega  \text{ is single-valued for }f \in
\Delta(\T)\,\}, 
\end{equation*}
then the decomposition \eqref{x-decomp-mesh} gives a decomposition
of the space $\P_r \Lambda^0(\T)$, i.e., the space of continuous
piecewise polynomials of degree $\le r$.

\section{Degrees of freedom and the dual decomposition}\label{sec:dof}

Although our main interest in this paper is obtaining direct sum
decompositions for polynomial differential forms that are analogous to the
Bernstein decomposition for ordinary polynomials, we include here a discussion
of another decomposition, referred to as the dual decomposition, for
completeness and as an illustration of the general theory developed in the
previous section.

Before we consider the case of differential forms, we review the
corresponding decomposition for polynomials.  For the construction of
finite element spaces based on the local space $\P_r(T)$, a basis for
the dual space $\P_r(T)^*$ is given, with each basis element
associated to a subsimplex of $T$.  This is referred to as a set of
\emph{degrees of freedom} for $\P_r(T)$.  The degrees of freedom then
determine the interelement continuity imposed on the finite element
space.  Indeed, in the classical approach of Ciarlet \cite{Ciarlet},
the degrees of freedom and their association to subsimplices is used
to define a finite element space.  For this purpose, what matters is
not the particular basis of $\P_r(T)^*$, but rather the decomposition
of this space into the spaces spanned by the basis elements associated
to each simplex.  For the standard Lagrange finite elements, this
geometric decomposition of the dual space is
\begin{equation}\label{dual-decomp}
\P_r(T)^* = \bigoplus_{f\in\Delta(T)} W_r(T,f),
\end{equation}
where
\[
W_r(g,f) := \{\,\psi\in \P_r(g)^*\, |\,
\psi(\omega)= \int_f (\tr_{g,f} \omega)\eta, \ 
  \eta\in \P_{r-\dim f-1}(f) \, \}.
\]
We note that for $\omega \in \P_r(h)$, $\tr_{h,f} \omega$ is uniquely
determined by $\bigoplus_{g\in\Delta(f)} W_r(h,g)$.

Consequently, if $f \subseteq h \in \Delta(T)$, we may define an
extension operator $F_{f,h} = F_{f,h}^r: \P_r(f) \to \P_r(h)$, determined
by the conditions:
\begin{gather*}
\int_g (\tr_{h,g} F_{f,h}\omega)\eta = \int_g (\tr_{f,g} \omega)\eta, 
\quad \eta\in\P_{r-\dim g-1}(g), \quad g\in\Delta(f),
\\
\psi(F_{f,h} \omega) =0, \quad
\psi \in W_r(h,g),\quad g\in\Delta(h),\quad g\nsubseteq f.
\end{gather*}

%Such a decomposition of the dual space implies a decomposition of
%the space $\P_r(T)$ itself.  For $f\in\Delta(T)$, we define
%\begin{gather}
%\label{Qr}
%\Q_r(T,f)= \{\, \omega\in\P_r(T)\,|\, \psi(\omega)=0 \text{ 
%for all $\psi\in\bigoplus_{g\nsubseteq f} W_r(T,g)$}\,\},
%\\
%\label{Qr0}
%\0\Q_r(T,f)= \{\, \omega\in\P_r(T)\,|\, \psi(\omega)=0 
%\text{ for all $\psi\in\bigoplus_{g\ne f} W_r(T,g)$}\,\},
%\end{gather}
%and note that $\0\Q_r(T,f)$, defined by \eqref{Qr0}, satisfies
%\[
%\0\Q_r(T,f) = \{ \omega \in \Q_r(T,f)\, |\, \tr_{T,f}\omega \in
%\0\P_r(f) \, \}.
%\]
To apply the theory developed in Section~\ref{cgd}, we need to check that the
extension operator is consistent, i.e., that it satisfies \eqref{H2}.  
For $f,g\subseteq h$, let $\omega \in \P_r(f)$, and set $\mu:=F_{f \cap g,g} \tr_{f, f \cap g}\omega \in \P_r(g)$,
$\nu:=\tr_{h,g} F_{f,h} \omega \in \P_r(g)$.  For any face 
$e\subseteq g\cap f$,
$\tr_{g,e}\mu=\tr_{f,e}\omega=\tr_{g,e}\nu$.  Therefore
$\psi(\mu)=\psi(\nu)$ for all $\psi\in\W_r(g,e)$ with $e\in\Delta(g)$ such
that
$e\subseteq f$.  Also, for $e\in\Delta(g)$ with $e\nsubseteq f$, it follows from the definition
of the extension that for all $\psi\in W_r(g,e)$, $\psi(\mu)=0=\psi(\nu)$.
Thus we have shown that the extension operators $F_{f,h}$ form a consistent family.
The decomposition
\begin{equation*}
\P_r(T) = \bigoplus_{f\in\Delta(T)} F_{f,T}[\0\P_r(f)],
\end{equation*}
corresponding to \eqref{x-decomp}, is now called the decomposition dual
to \eqref{dual-decomp}. Furthermore, from Theorem~\ref{basic-decomp}
we obtain a 
corresponding direct sum decomposition for the assembled space $\P_r(\T)
= \P_r\Lambda^0(\T)$ of the form \eqref{x-decomp-mesh}.

%Assumption~\ref{A1} follows from the fact that if $f \subseteq g$, then
%\begin{equation*}
%\bigoplus_{h\nsubseteq g} W_r(S,h) \subseteq 
%\bigoplus_{h\nsubseteq f} W_r(S,h).
%\end{equation*}

%The fact that properties Assumption~\ref{A0}--\ref{A3} are true will follow from the discussion
%given in Section~\ref{sec:dof} below.

In the remainder of this section, we present analogous results for the spaces
$\P_r\Lambda^k(T)$ and $\P^-_r\Lambda^k(T)$. This will be based on the
following decompositions of the dual spaces 
$\P_r\Lambda^k(T)^*$ and $\P^-_r\Lambda^k(T)^*$,
established in \cite[\S 4, Theorems 4.10 and 4.14]{acta}.

\begin{thm}\label{stardecomp}
1. For each $f \in \Delta(T)$ define
\begin{equation*}
W_r^k(T,f) := \biggl\{\,\psi\in \P_r\Lambda^k(T)^*\,\big|\,
\  \psi(\omega)= \int_f \tr_{T,f}\omega\wedge\eta \text{ for some $\eta\in
   \P_{r+k-\dim f}^-\Lambda^{\dim f-k}(f)$}\,\biggr\}.
\end{equation*}
Then the obvious mapping $\P_{r+k-\dim f}^-\Lambda^{\dim f-k}(f)\to W_r^k(T,f)$
is an isomorphism, and
$$
\P_r\Lambda^k(T)^* = \bigoplus_{f\in\Delta(T)} W_r^k(T,f).
$$

2. For each $f \in \Delta(T)$ define
\begin{equation*}
W_r^{k-}(T,f) := \biggl\{\,\psi\in \P^-_r\Lambda^k(T)^*\,\big|\,
\ \psi(\omega)= \int_f \tr_{T,f}\omega\wedge\eta \text{ for some $\eta\in
   \P_{r+k-\dim f-1}\Lambda^{\dim f-k}(f)$}\,\biggr\}.
\end{equation*}
Then the obvious mapping $\P_{r+k-\dim f-1}\Lambda^{\dim f-k}(f)\to W_r^{k-}(T,f)$
is an isomorphism, and
$$
\P^-_r\Lambda^k(T)^* = \bigoplus_{f\in\Delta(T)} W_r^{k-}(T,f).
$$
\end{thm}
Note that as in the polynomial case, if $\omega \in \P_r\Lambda^k(T)$, then
$\tr_{T,f} \omega$ is determined by the degrees of freedom in $W_r^{k}(T,g)$
for $g \in \Delta(f)$.  In particular, if $\omega \in \P_r\Lambda^k(T)$ such
that all the degrees of freedom associated to the subsimplices of $T$ with
dimension less than or equal to $n-1$ vanish, then $\omega \in
\0\P_r\Lambda^k(T)$.  The corresponding property holds for the
spaces $\P^-_r\Lambda^k(T)$ as well.

An immediate consequence of this theorem are the following isomorphisms,
that will be used in the following section.
\begin{cor}\label{iso-spaces}
\[
\0\P_r\Lambda^k(T)^* \cong \P^-_{r+k-n}\Lambda^{n-k}(T) \quad
\text{and} \quad \0\P^-_r\Lambda^k(T)^* \cong \P_{r+k-n-1}\Lambda^{n-k}(T).
\]
\end{cor}

As in the case of 0-forms, if $f \subset h \in
\Delta(T)$, we  define an extension operator $F_{f,h}^{k,r} : \P_r
\Lambda^k(f) \to \P_r\Lambda^k (h)$, determined by the conditions:
\begin{gather*}
\int_g(\tr_{h,g}F_{f,h}^{k,r}\omega)\eta =
\int_g(\tr_{f,g}\omega)\eta,\quad \eta\in\P_{r-\dim g-1}(g), \quad g\in \Delta(f),
\\
\psi(F_{f,h}^{k,r} \omega) =0,\quad 
\psi \in  W_r^k(h,g),\quad g\in\Delta(h), \quad g\nsubseteq f.
\end{gather*}
We may similarly define an extension operator $F_{f,T}^{k,r,-} : \P_r^- \Lambda^k(f) \to
\P_r^- \Lambda^k (h)$.  The verification of the consistency of these families of extension
operators is essentially the same as for the space $\P_r(T)$ given above,
and so we do not repeat the proof.

%If $\T$ is a simplicial triangulation of a domain $\Omega \in
%\R^n$, then the space $\P_r\Lambda^k(\T)$ can alternatively be
%characterized as
%\begin{multline*}
%\P_r\Lambda^k(\T) =\{\,\omega\in L^2\Lambda^k(\Omega)\,:\, 
%\omega|_T\in\P_r\Lambda^k(T) \quad \forall T\in\T,
%\\*[-3pt] 
% \phi(\omega) \, \text{is singled-valued for }f \in
%\Delta_j(\T),\, k \le j \le n-1, \, \phi \in W_r^{k}(f) \, \},
%\end{multline*}
%where 
%\[
%W_r^{k}(f) = \bigcap_{\substack{T \in \T\\f \in
%    \Delta(T)}}W_r^{k}(T,f).
%\]
%A completely analogous characterization can be given for the spaces
%$\P^-_r\Lambda^k(T)$.
%We have therefore seen that the

\section{Barycentric spanning sets}\label{sec:span-basis}
Let $T=[x_0,\ldots ,x_n] \subset \R^n$ be a nondegenerate $n$-simplex. 
The Bernstein basis
described in Section~\ref{Bernstein} above is given in terms of the
barycentric coordinates $\{\lambda_i\}_{i=0}^n \subset \P_1(T)$. The
main purpose of this paper is to construct the generalization of the Bernstein basis
for 
the polynomial spaces $\P_r\Lambda^k(T)$
and $\P^-_r\Lambda^k(T)$. In the present section, we will give spanning
sets and bases for these spaces and for the corresponding spaces with vanishing
trace expressed in barycentric coordinates.  Note that the bases given
in this section depend on the ordering of the vertices.  These are
\emph{not} the bases we suggest for computation.

For convenience we summarize the results of the section in the following theorem,
referring not only to the $n$-dimensional simplex $T$, but, more generally, to
any subsimplex $f$ of $T$.
Here, we use the notation
\begin{equation}\label{dlsigma}
 d\lambda^f_\sigma = d\lambda^f_{\sigma(1)} \wedge \cdots \wedge 
d\lambda^f_{\sigma(k)} \in \Alt^k T_f
\end{equation}
for $f \in \Delta(T)$, $\sigma \in \Sigma(1:k, 0:n)$ with
$\llbracket\sigma\rrbracket \subseteq \I(f)$,
 and $\phi^f_\sigma$ for the Whitney form defined in \eqref{wf}.

\begin{thm}\label{summary}
Let $f\in\Delta(T)$.

1. Spanning set and basis for $\P_r\Lambda^k(f)$.  The set
\begin{equation*}
\{\,(\lambda^f)^\alpha d\lambda^f_\sigma\,|\,\alpha\in\N_0^{0:n},|\alpha|=r,
\sigma\in\Sigma(1:k,0:n),\ \llbracket\alpha, \sigma\rrbracket\subseteq \I(f)\,\}
\end{equation*}
is a spanning set for $\P_r\Lambda^k(f)$, and
\begin{equation*}
\{\,(\lambda^f)^\alpha d\lambda^f_\sigma\,|\,\alpha\in\N_0^{0:n},|\alpha|=r,
\sigma\in\Sigma(1:k,0:n),\ \llbracket\alpha, \sigma\rrbracket\subseteq \I(f),
\ \min\llbracket\sigma\rrbracket>\min\I(f) \,\}
\end{equation*}
is a basis.

2. Spanning set and basis for $\0\P_r\Lambda^k(f)$.  The set
\begin{equation*}
\{\,(\lambda^f)^\alpha d\lambda^f_\sigma\,|\,\alpha\in\N_0^{0:n},|\alpha|=r,
\sigma\in\Sigma(1:k,0:n),\ \llbracket\alpha, \sigma\rrbracket=\I(f)\,\}
\end{equation*}
is a spanning set for $\0\P_r\Lambda^k(f)$, and
\begin{equation*}
\{\,(\lambda^f)^\alpha d\lambda^f_\sigma\,|\,\alpha\in\N_0^{0:n},|\alpha|=r,
\sigma\in\Sigma(1:k,0:n),\ \llbracket\alpha, \sigma\rrbracket=\I(f)
\ \alpha_i=0\text{ if } i<\min[\I(f)\setminus\llbracket\sigma\rrbracket] \,\}
\end{equation*}
is a basis.

3. Spanning set and basis for $\P_r^-\Lambda^k(f)$.  The set
\begin{equation*}
\{\,(\lambda^f)^\alpha \phi^f_\sigma\,|\,\alpha\in\N_0^{0:n},|\alpha|=r-1,
\sigma\in\Sigma(0:k,0:n),\ \llbracket\alpha, \sigma\rrbracket\subseteq \I(f)\,\}
\end{equation*}
is a spanning set for $\P^-_r\Lambda^k(f)$, and
\begin{equation}\label{prmfbas}
\{\,(\lambda^f)^\alpha \phi^f_\sigma\,|\,\alpha\in\N_0^{0:n},|\alpha|=r-1,
\sigma\in\Sigma(0:k,0:n),\ \llbracket\alpha, \sigma\rrbracket\subseteq \I(f),
\ \alpha_i=0\text{ if } i<\min\llbracket\sigma\rrbracket \,\}
\end{equation}
is a basis.

4. Spanning set and basis for $\0\P_r^-\Lambda^k(f)$.  The set
\begin{equation*}
\{\,(\lambda^f)^\alpha \phi^f_\sigma\,|\,\alpha\in\N_0^{0:n},|\alpha|=r-1,
\sigma\in\Sigma(0:k,0:n),\ \llbracket\alpha, \sigma\rrbracket= \I(f)\,\}
\end{equation*}
is a spanning set for $\0\P^-_r\Lambda^k(f)$, and
\begin{equation*}
\{\,(\lambda^f)^\alpha \phi^f_\sigma\,|\,\alpha\in\N_0^{0:n},|\alpha|=r-1,
\sigma\in\Sigma(0:k,0:n),\ \llbracket\alpha, \sigma\rrbracket=\I(f),
\ \alpha_i=0\text{ if } i<\min\llbracket\sigma\rrbracket \,\}
\end{equation*}
is a basis.
\end{thm}

\subsection{Barycentric spanning set and basis for $\P_r\Lambda^k(T)$}
\label{subsec:spanpr}
Observe that $d\lambda_i \in \Alt^1 \R^n$. Furthermore,
$d\lambda_i(x_j - y) = \delta_{ij}$ for any $y$ in the subsimplex
opposite $x_i$. In particular, $\tr_{T,f} d\lambda_i = 0$ for any
subsimplex $f \in \Delta(T)$ with $x_i \notin f$ or equivalently $i
\in \I(f^*)$.
Furthermore, $\{d\lambda_i\}_{i=0}^n$ is a spanning set
for $\Alt^1 \R^n$, and any subset of $n$ elements is a basis.
Therefore, writing $d\lambda_\sigma$ for $d\lambda^T_\sigma$, the set
\[
\{\, d\lambda_\sigma\,|\, \sigma \in \Sigma(1:k,0:n)\, \}
\]
is a spanning set for $\Alt^k \R^n$, and the set
\[
\{\, d\lambda_\sigma \,|\, \sigma \in \Sigma(1:k,1:n)\, \}
\]
is a basis. The forms $d\lambda_\sigma \in \Alt^k
\R^n$
have the property that for any $f \in \Delta(T)$,
\[
\tr_{T,f} d\lambda_\sigma = 0 \quad \text{if and only if } \llbracket
\sigma \rrbracket \cap \I(f^*) \neq \emptyset.
\]
More generally, for polynomial forms of the form $\lambda^\alpha d\lambda_\sigma
\in \P_r\Lambda^k(T)$, we observe that 
\[
\tr_{T,f} (\lambda^\alpha d\lambda_\sigma) = 0 \quad \text{if and only if } \llbracket
\alpha,\sigma \rrbracket \cap \I(f^*) \neq \emptyset.
\]
In particular, $\lambda^\alpha d\lambda_\sigma \in \0\P_r\Lambda^k(T)$
if and only if $\llbracket
\alpha,\sigma \rrbracket = \{0,\ldots ,n\}$.

Taking the tensor product of the Bernstein basis for $\P_r(T)$, given by
\eqref{span1}, with the spanning set and basis given above for $\Alt^k \R^n$, we get
that
\begin{prop}\label{prspan}
The set
\begin{equation*}
\{\,\lambda^\alpha d\lambda_\sigma\,|\,\alpha\in\N_0^{0:n},|\alpha|=r,
\sigma\in\Sigma(1:k,0:n)\,\}
\end{equation*}
is a spanning set for $\P_r\Lambda^k(T)$, and
\begin{equation*}
\{\,\lambda^\alpha d\lambda_\sigma\,|\,\alpha\in\N_0^{0:n},|\alpha|=r,
\sigma\in\Sigma(1:k,1:n)\,\}
\end{equation*}
is a basis.
\end{prop}
Restricting to a face $f \in \Delta(T)$, we obtain the spanning set and basis
for $\P_r \Lambda^k(f)$ given in the first part of Theorem~\ref{summary}.

\subsection{Barycentric spanning set and basis for $\P^-_r\Lambda^k(T)$}
\label{subsec:spanprminus}
For $f \in \Delta(T)$ and $\sigma\in\Sigma(0:k,0:n)$ with
$\llbracket\sigma\rrbracket \subseteq \I(f)$, define
the associated Whitney form by
\begin{equation}\label{wf}
\phi^f_\sigma = \sum_{i=0}^k (-1)^i\lambda^f_{\sigma(i)}\,
d\lambda^f_{\sigma(0)}\wedge\dots \wedge
\widehat{d\lambda^f_{\sigma(i)}}\wedge\dots \wedge d\lambda^f_{\sigma(k)}.
\end{equation}
Just as we usually write $\lambda_i$ rather than $\lambda^T_i$ when the simplex is clear from context, we will usually write $\phi_\sigma$ instead of $\phi^T_\sigma$.
We note that if $k=0$, so that the associated subsimplex $f_\sigma$
consists of a single point $x_i$, then $\phi_\sigma = \lambda_i$.
It is evident that the Whitney forms belong to
$\P_1\Lambda^k(T)$. In fact, they belong to
$\P^-_1\Lambda^k(T)$.
This is a direct consequence of the identity
\begin{equation}\label{phi-id-1}
\kappa d\lambda_\sigma = \phi_\sigma - \phi_\sigma(0),
\end{equation}
which can be easily established by induction on $k$,
using the Leibniz rule \eqref{kappa-leibniz}.
In fact, the set
\[
\{\, \phi_\sigma \,|\, \sigma \in \Sigma(0:k,0:n) \, \}
\]
is a basis for $\P_1^-\Lambda^k(T)$. Furthermore,
$\tr_{T,f} \phi_\sigma = d\lambda_{\sigma(1)} \wedge \cdots \wedge
d\lambda_{\sigma(k)}$ 
is a nonvanishing constant $k$-form on $f=f_\sigma$, while $\tr_{T,f} \phi_\sigma = 0$
for $f \in \Delta_k(T)$, $f \neq f_\sigma$. Therefore,
we refer to $\phi_\sigma$ as \emph{the Whitney form associated to the face
$f_\sigma$}.

For $\sigma\in\Sigma(0:k,0:n)$ and $0 \le j \le k$,
we let $\phi_{\sigma\hat\jmath}$ be the Whitney form corresponding to the
subsimplex of $f_\sigma$ obtained by removing the vertex $\sigma(j)$.
Hence,
\begin{align*}
\phi_{\sigma\hat\jmath} & = \sum_{i=0}^{j-1} (-1)^i\lambda_{\sigma(i)}
d\lambda_{\sigma(0)}\wedge\dots \wedge
\widehat{d\lambda_{\sigma(i)}}\wedge\dots \wedge 
\widehat{d\lambda_{\sigma(j)}}\wedge\dots
d\lambda_{\sigma(k)}\\
&-\sum_{i=j+1}^{k} (-1)^i\lambda_{\sigma(i)}
d\lambda_{\sigma(0)}\wedge\dots \wedge
\widehat{d\lambda_{\sigma(j)}}\wedge\dots \wedge 
\widehat{d\lambda_{\sigma(i)}}\wedge\dots
d\lambda_{\sigma(k)}.
\end{align*}
From this expression, we easily obtain the identity
\begin{equation}\label{whitid}
\sum_{j=0}^k(-1)^j\lambda_{\sigma(j)}\phi_{\sigma\hat\jmath} =0,
\quad \sigma\in\Sigma(0:k,0:n).
\end{equation}
Correspondingly, for $j \notin \llbracket \sigma \rrbracket$,
we define
\[
\phi_{j\sigma} = \lambda_j d\lambda_\sigma - d\lambda_j \wedge \phi_\sigma.
\]
Thus, modulo a possible factor of $-1$,
$\phi_{j\sigma}$ is the Whitney form associated to the simplex
$[x_j,f_\sigma]$.
For these functions, we obtain 
\begin{equation}\label{phi-id-2}
\begin{aligned}
 \sum_{j \notin \llbracket \sigma \rrbracket}\phi_{j\sigma}
&= (\sum_{j \notin \llbracket \sigma \rrbracket}
\lambda_j) d\lambda_\sigma - (\sum_{j \notin \llbracket \sigma \rrbracket}
d\lambda_j) \wedge \phi_\sigma
\\
&= (\sum_{j \notin \llbracket \sigma \rrbracket}
\lambda_j) d\lambda_\sigma
+ (\sum_{j \in \llbracket \sigma \rrbracket}
d\lambda_j) \wedge \phi_\sigma
= (\sum_{j=0}^n \lambda_j)d\lambda_\sigma = d\lambda_\sigma.
\end{aligned}
\end{equation}
Now consider functions of the form $\lambda^\alpha \phi_\sigma$, where
$\alpha\in\N_0^{0:n}$, $|\alpha|=r-1$,
$\sigma\in\Sigma(0:k,0:n)$. It follows from the relation
\eqref{closure} that these 
functions belong to $\P^-_r\Lambda^k(T)$.  In fact, they span.
From the identity \eqref{whitid} we know that these
forms are not, in general, linearly independent. The following
lemma, cf.~\cite[Lemma 4.2]{acta}, enables us to extract a basis.
\begin{lem}\label{linind}
Let $x$ be a vertex of $T$.  Then the Whitney forms corresponding
to the $k$-subsimplices that contain $x$ are linearly independent
over the ring of polynomials $\P(T)$.
\end{lem}
Using these results, we are able to prove:
\begin{prop}\label{prminusspan}
The set
\begin{equation}\label{spanning-P-}
\{\,\lambda^\alpha \phi_\sigma\,|\,\alpha\in\N_0^{0:n},|\alpha|=r-1,
\sigma\in\Sigma(0:k,0:n)\,\}.
\end{equation}
is a spanning set for $\P^-_r\Lambda^k(T)$, and
\begin{equation}
\label{pr--basis}
\{\,\lambda^\alpha \phi_\sigma\,|\,\alpha\in\N_0^{0:n},\ |\alpha|=r-1,\
\sigma\in\Sigma(0:k,0:n),\ \alpha_i=0\text{ if } i<\min\llbracket\sigma\rrbracket \,\}
\end{equation}
is a basis.
\end{prop}
\begin{proof}
Let $\alpha\in\N_0^{0:n}$, with $|\alpha|=r-1$, and $\rho \in
\Sigma(0:k-1,0:n)$.
The identity \eqref{phi-id-2} implies that
\[
\lambda^\alpha d\lambda_\rho = \sum_{j \notin \llbracket \rho
  \rrbracket}
\lambda^\alpha \phi_{j\rho},
\]
and hence all forms in $\P_{r-1}\Lambda^k(T)$ are in the span of the set
given by \eqref{spanning-P-}. Furthermore, if 
$\sigma \in \Sigma(0:k,0:n)$, we obtain from
\eqref{phi-id-1} that
\[
\kappa(\lambda^\alpha d\lambda_\sigma) + \lambda^\alpha \phi_\sigma(0)
= \lambda^\alpha \phi_\sigma,
\]
and therefore all of $\kappa[\P_{r-1}\Lambda^{k+1}(T)]$ is also in the
span.  By the definition of the space $\P^-_r\Lambda^k(T)$, it follows
that \eqref{spanning-P-} is a spanning set.  To show that
\eqref{pr--basis} is a basis, we use the identity \eqref{whitid} to see
that any form given in the span of \eqref{spanning-P-} is in
the span of the forms in \eqref{pr--basis}.  Then we use
Lemma~\ref{linind}, combined with a simple inductive argument, to show
that the elements of the asserted basis are linearly independent.  For
details, see the proof of Theorem 4.4 of \cite{acta}.
\end{proof}

Restricting to a face $f \in \Delta(T)$, we obtain the spanning set and basis
for $\P_r^- \Lambda^k(f)$ given in the third part of Theorem~\ref{summary}.

\subsection{Spaces of vanishing trace}
In this subsection, we will derive spanning sets and bases for
the corresponding spaces of zero trace.  This will be based on
the results obtained above and Corollary~\ref{iso-spaces}, which
leads to the dimension of these spaces.  We first
characterize the space $\0 \P^-_r\Lambda^k(T)$.
\begin{prop} \label{0spacespan-}
The set
\begin{equation*}
\{\,\lambda^\alpha \phi_\sigma\,|\,\alpha\in\N_0^{0:n},|\alpha|=r-1,
\sigma\in\Sigma(0:k,0:n),
\llbracket\alpha, \sigma\rrbracket=\{0,\ldots,n\}\,\}
\end{equation*}
is a spanning set for $\0\P^-_r\Lambda^k(T)$ and 
\begin{equation*}
\{\,\lambda^\alpha \phi_\sigma\,|\,\alpha\in\N_0^{0:n},\ |\alpha|=r-1,\
\sigma\in\Sigma(0:k,0:n),\ 
\llbracket\alpha, \sigma\rrbracket=\{0,\ldots,n\},
\ \alpha_i=0\text{ if } i<\min\llbracket\sigma\rrbracket\,\}
\end{equation*}
is a basis.
\end{prop}
\begin{proof}
Since $\llbracket\alpha,\sigma \rrbracket = \{0,\ldots,n\}$, each of the forms
$\lambda^\alpha \phi_\sigma$ is contained in $\0 \P^-_r\Lambda^k(T)$.
Moreover, the condition $\alpha_i=0$ if $i<\min\llbracket\sigma\rrbracket$ reduces to $\sigma(0)=0$ in this case.
Lemma~\ref{linind} implies that
the forms $\lambda^{\alpha} \phi_\sigma$ for which
$\sigma(0)= 0$ are linearly independent.  The cardinality of this
set
is equal to $\binom{n}{k}\dim \P_{r-n+k-1}$ which is equal to 
$\dim \0\P^-_r\Lambda^k(T)$ by Corollary~\ref{iso-spaces}.
This completes the proof.
\end{proof}

Restricting to a face $f \in \Delta(T)$, we obtain the spanning set and basis
for $\0\P_r^- \Lambda^k(f)$ given in the fourth part of Theorem~\ref{summary}.

Finally, we obtain a characterization of the space $\0 \P_r\Lambda^k(T)$.
\begin{prop} \label{0spacespan}
The set 
\begin{equation*}
\{\,\lambda^\alpha d\lambda_\sigma\,|\,\alpha\in\N_0^{0:n},\ |\alpha|=r,
\ \sigma\in\Sigma(1:k,0:n), \ 
\llbracket\alpha,\sigma \rrbracket = \{0,\ldots,n\}\,\}
\end{equation*}
is a spanning set for $\0 \P_r\Lambda^k(T)$, and
\begin{equation}
\label{pr0-basis} 
\{\,\lambda^\alpha d\lambda_\sigma\,|\,
\alpha\in\N_0^{0:n},\ |\alpha|=r,
\ \sigma\in\Sigma(1:k,0:n), \
\llbracket\alpha,\sigma \rrbracket = \{0,\ldots,n\},
\
\alpha_i = 0 \text{ if } i < \min\llbracket\sigma^*\rrbracket \,\}
\end{equation}
is a basis.
\end{prop}

\begin{proof}
Since $\llbracket\alpha,\sigma \rrbracket = \{0,\ldots,n\}$, each of
the forms $\lambda^\alpha d\lambda_\sigma$ is contained in $\0
\P_r\Lambda^k(T)$. Furthermore, we have seen in
Corollary~\ref{iso-spaces}, that $\0\P_r\Lambda^k(T)
\cong\P^-_{r+k-n}\Lambda^{n-k}(T)$, whence, $\dim \0\P_r\Lambda^k(T) =
\binom{r-1}{n-k} \binom{r+k}{r}$.  On the other hand, the cardinality
of the set given by \eqref{pr0-basis} can be computed as $\sum_j
A_j\cdot B_j$, where $A_j$ is the number of elements $\sigma \in
\Sigma(1:k,0:n)$ with $\min\llbracket\sigma^*\rrbracket = j$, and for
each fixed such $\sigma$, $B_j$ is the number of multi-indices
$\alpha$ satisfying the conditions of \eqref{pr0-basis}, namely
\[
A_j = \binom{n-j}{k-j} \quad \text{and} \quad B_j =
\binom{r+k-j-1}{n-j}.
\]
Hence, the cardinality of the set is given by 
\begin{align*}
\sum_{j=0}^k \binom{n-j}{k-j}\binom{r+k-j-1}{n-j}
&= \binom{r-1}{n-k} \sum_{j=0}^k \binom{r+k-j-1}{r-1}\\
&=\binom{r-1}{n-k}\binom{r+k}{r}.
\end{align*}
Here the first identity follows from a binomial identity of the
form
\[
\binom{a}{b} \binom{b}{c} = \binom{a}{c} \binom{a-c}{b-c},
\]
while the second is a standard summation formula.
Hence, the cardinality of the set given by
\eqref{pr0-basis} is equal to the dimension of $\0\P_r\Lambda^k(T)$.
To complete the proof, we show that the elements of the set \eqref{pr0-basis} are linearly independent.   Denote the index set by
\begin{equation*}
S:=\{\,(\alpha,\sigma)\in\N^{0:n}_0\x \Sigma(1:k,0:n)\,|\,
|\alpha|=r,\ \llbracket\alpha,\sigma \rrbracket = \{0,\ldots,n\},
\ \alpha_i = 0 \text{ if } i < \min\llbracket\sigma^*\rrbracket \,\},
\end{equation*}
so we must show that if
\begin{equation}\label{lin-indep?}
\sum_{(\alpha,\sigma)\in S} c_{\alpha\sigma} \lambda^\alpha d\lambda_\sigma
= 0,
\end{equation}
for some real coefficients $c_{\alpha\sigma}$, then all the coefficients vanish.  Since the Bernstein monomials $\lambda^\alpha$
are linearly independent, \eqref{lin-indep?} implies that
for each $\alpha\in\N^{0:n}_0$ with $|\alpha|=r$,
\begin{equation}\label{lin-ind-alpha}
\sum_{\{\,\sigma\,|\,(\alpha,\sigma)\in S\,\}} c_{\alpha\sigma}  d\lambda_\sigma =0.
\end{equation}
First consider a multi-index $\alpha$ with $\alpha_0>0$.  Then
the definition of the index set $S$ implies that $\min\llbracket\sigma^*\rrbracket=0$
for all the summands in \eqref{lin-ind-alpha}.  Since the corresponding $d\lambda_\sigma$ are linearly independent, we
conclude that all the $c_{\alpha\sigma}$ vanish when $\alpha_0>0$.
Next consider $\alpha$ with $\alpha_0=0$ but $\alpha_1>0$.
If $(\alpha,\sigma)\in S$, then $\min\llbracket\sigma^*\rrbracket=1$, and again we
conclude that $c_{\alpha\sigma}=0$.  Continuing in this way
we find that all the $c_{\alpha\sigma}$ vanish, completing the proof.
\end{proof}

Restricting to a face $f \in \Delta(T)$, we obtain 
the spanning set and basis
for $\0\P_r \Lambda^k(f)$ given in the second part of Theorem~\ref{summary}.

\section{A geometric decomposition of $\P^-_r\Lambda^k(T)$}\label{sec:decomp-}

In this section, we will apply the theory developed in Section~\ref{cgd} with
$X(T) = \P^-_r\Lambda^k(T)$ to obtain a geometric
decomposition of $\P_r\Lambda^k(\T)$ into subspaces $E_f[\0\P^-_r\Lambda^k(f)]$,
where $E_f$ is the global extension operator constructed as in Section~\ref{cgd} from a  consistent family
of easily computable extension operators.  The resulting decomposition
reduces to the Bernstein decomposition \eqref{B-decomp} in the case
$k=0$.

We first note that if $T,T'\in\T$ with $f\in\Delta(T)\cap\Delta(T')$
then $\tr_{T,f}\P_r^-\Lambda^k(T)=\tr_{T',f}\P_r^-\Lambda^k(T') =
\P^-_r\Lambda^k(f)$.  Hence, the assumption \eqref{cassump} holds.
Furthermore, for $f, g \in \Delta(\T)$ with $f \subseteq g$, we define
$E =E_{f,g}^{k,r,-}: \P^-_r\Lambda^k(f) \to \P^-_r\Lambda^k(g)$ as the
barycentric extension:
\begin{equation}
\label{bary-ext-}
(\lambda^f)^\alpha \phi^f_\sigma \mapsto (\lambda^g)^\alpha
\phi^g_\sigma, \qquad \llbracket\alpha,\sigma\rrbracket \subseteq \I(f).
\end{equation}
This generalizes to $k$-forms, the barycentric extension operator $E^r_{f,T}$
on $\P_r$, introduced in Section~\ref{Bernstein}. 
Since the forms $(\lambda^f)^\alpha \phi^f_\sigma$ are not linearly
independent, it is not clear that \eqref{bary-ext-} well-defines $E$.
We show this  in the following theorem.
\begin{thm}\label{P--spaces-tr}
There is a unique mapping $E=E_{f,g}^{k,r,-}$ from $\P^-_r\Lambda^k(f)$ to $\P^-_r\Lambda^k(g)$
satisfying \eqref{bary-ext-}.
\end{thm}
\begin{proof}
We first recall from part 3 of Theorem~\ref{summary} that the set
\begin{equation*}
 \{\,(\lambda^f)^\alpha
\phi^f_\sigma\,|\,\alpha\in\N_0^{0:n},|\alpha|=r-1, \sigma\in\Sigma(0:k,0:n),
\ \llbracket\alpha,\sigma\rrbracket \subseteq \I(f), 
\ \alpha_i=0\text{ if } i<\min\llbracket\sigma\rrbracket\,\}
\end{equation*}
is a basis for $\P^-_r\Lambda^k(f)$. Hence, we can uniquely define an extension $E$ by 
\eqref{bary-ext-}, if we restrict to these basis functions.
 We now show that \eqref{bary-ext-} holds for all $(\lambda^f)^\alpha
\phi^f_\sigma$
with $\llbracket\alpha,\sigma\rrbracket \subseteq \I(f)$.
%% To see this, we note that the identity \eqref{whitid} 
%% implies that
%% any such $(\lambda^f)^\alpha \phi^f_{\sigma}$ belongs to 
%% the span of the basis functions. In fact, if
%% $(\lambda^f)^\alpha \phi^f_{\sigma} = \sum c_{\beta, \rho}
%% (\lambda^f)^\beta \phi^f_{\rho}$ on $f$, where the sum runs over the
%% basis functions, then $(\lambda^g)^\alpha \phi^g_{\sigma} = \sum c_{\beta, \rho}
%% (\lambda^g)^\beta \phi^g_{\rho}$ on $g$.
%% Hence,
%% \begin{equation*}
%% E[(\lambda^f)^{\alpha} \phi^f_{\sigma}] = \sum c_{\beta, \rho}
%% E[(\lambda^f)^\beta \phi^f_{\rho}]
%% = \sum c_{\beta, \rho} (\lambda^g)^\beta \phi^g_\rho
%% = (\lambda^g)^\alpha \phi^g_\sigma.
%% \end{equation*}
To see this, we use the identity \eqref{whitid}. Consider forms
$(\lambda^f)^\alpha \phi^f_{\sigma}$ which do not belong to the given
basis,
i.e., $s := \min\llbracket \alpha \rrbracket <\min
\llbracket \sigma \rrbracket$. Write $\lambda^\alpha = 
\lambda^\beta\lambda_s$ and let $\rho \in \Sigma(0:k+1,0:n)$ be determined by 
$\llbracket \rho \rrbracket = \{s\}\cup \llbracket \sigma \rrbracket$.
Then, by \eqref{whitid},
\[
(\lambda)^\alpha \phi_{\sigma} = 
\sum_{j=1}^{k+1}(-1)^{j-1}\lambda^\beta\lambda_{\rho(j)}\phi_{\rho\hat\jmath},
\]
and so
\[
(\lambda^f)^\alpha \phi^f_{\sigma} = 
\sum_{j=1}^{k+1}(-1)^{j-1}(\lambda^f)^\beta\lambda^f_{\rho(j)}\phi^f_{\rho\hat\jmath}, \qquad
(\lambda^g)^\alpha \phi^g_{\sigma} = 
\sum_{j=1}^{k+1}(-1)^{j-1}(\lambda^g)^\beta\lambda^g_{\rho(j)}\phi^g_{\rho\hat\jmath}.
\]
Hence
\begin{equation*}
E[(\lambda^f)^{\alpha} \phi^f_{\sigma}] =
\sum_{j=1}^{k+1}(-1)^{j-1}E[(\lambda^f)^\beta
\lambda^f_{\rho(j)}\phi^f_{\rho\hat\jmath}]
=
\sum_{j=1}^{k+1}(-1)^{j-1}(\lambda^g)^\beta
\lambda^g_{\rho(j)}\phi^g_{\rho\hat\jmath}
= (\lambda^g)^\alpha \phi^g_\sigma,
\end{equation*}
and the proof is completed.
\end{proof}

\begin{thm}
\label{E-consistent}
The family of extension operators 
$E$ is consistent, i.e., for all $f,g,h\in\Delta(\T)$ with 
$f,g \subseteq h$, and all 
$\omega \in \P^-_r\Lambda^k(f)$,
\begin{equation*}
\tr_{h,g}E_{f,h} \omega = E_{f \cap g, g} \tr_{f,f \cap g}\omega.
\end{equation*}
\end{thm}
\begin{proof}
It is enough to establish this result  for $\omega = (\lambda^f)^{\alpha}
\phi^f_{\sigma}$, with $\llbracket\alpha,\sigma\rrbracket\subseteq\I(f)$, since such $\omega$
span $\P_r^-\Lambda^k(f)$.
Now for such pairs $(\alpha, \sigma)$, 
$E_{f,h}[(\lambda^f)^{\alpha} \phi^f_{\sigma}] = (\lambda^h)^{\alpha}
\phi^h_{\sigma}$ and then
\begin{equation*}
\tr_{h,g}E_{f,h}[(\lambda^f)^{\alpha} \phi^f_{\sigma}] = \begin{cases}
(\lambda^g)^{\alpha} \phi^g_{\sigma}, & \text{if } 
\llbracket\alpha,\sigma\rrbracket\subseteq\I(f \cap g),
\\ 0, & \text{otherwise}. \end{cases}
\end{equation*}
On the other hand,
\begin{equation*}
\tr_{f,f \cap g} (\lambda^f)^{\alpha} \phi^f_{\sigma}
= \begin{cases}
(\lambda^{f \cap g})^{\alpha} (\phi^{f \cap g})_{\sigma}, & \text{if }
\llbracket\alpha,\sigma\rrbracket\subseteq\I(f \cap g), \\
0, & \text{otherwise}, \end{cases}
\end{equation*}
and hence
\begin{equation*}
E_{f \cap g, g} \tr_{f,f \cap g} (\lambda^f)^{\alpha} \phi^f_{\sigma}
= \begin{cases}
(\lambda^{g})^{\alpha} \phi^{g}_{\sigma}, & \text{if }
\llbracket\alpha,\sigma\rrbracket\subseteq\I(f \cap g), \\
0, & \text{otherwise}. \end{cases}
\end{equation*}
\end{proof}

From Theorem~\ref{basic-decomp}, we obtain the desired geometric decomposition of $\P^-_r\Lambda^k(\T)$.
\begin{thm}
\label{direct-sum-decomp-}
\begin{equation*}
\P^-_r\Lambda^k(\T) = \bigoplus_{\substack{f\in\Delta(\T)
\\ \dim f\ge k}} E_f[\0\P^-_r\Lambda^k(f)].
\end{equation*}
where $E_f:\0\P^-_r\Lambda^k(f)\to \P^-_r\Lambda^k(\T)$ denotes the global extension
operator determined by the family $E_{f,g}^{k,r,-}$.
\end{thm}

The final part of Theorem~\ref{summary} furnishes an explicit spanning set and basis for $\0\P^-_r\Lambda^k(f)$, and
so this theorem gives an explicit spanning set and basis for $\P^-_r\Lambda^k(\T)$. We discuss these explicit representations further
in Section~\ref{sec:bases}.

We now turn to a geometric characterization of the extension operator
$E:\P^-_r\Lambda^k(f)\to\P^-_r\Lambda^k(T,f)$.  To this end, we say that
a smooth $k$-form $\omega\in \Lambda^k(T)$ vanishes to order $r$ at a point
$x$ if the function $x\mapsto\omega_x(v_1, \ldots ,v_k)$ vanishes to order $r$
at $x$ for all $v_1,\ldots ,v_k \in \R^n$, and that it vanishes to order $r$
on a set $g$ if it vanishes to order $r$ at each point of the set.
Note that the extension operator $E=E_{f,T}^{k,r,-}$ has the property
that for any $\mu\in\P^-_r\Lambda^k(f)$,
$E_{f,T}^{k,r,-} \mu$ vanishes to order $r$ on $f^*$.  In fact,
if we set
\begin{equation*}
\P^-_r\Lambda^k(T,f) = \{\, \omega \in \P^-_r\Lambda^k(T) \,|\,
\omega \text{ vanishes to order } r \text{ on $f^*$} \, \},
\end{equation*}
we can prove 
\begin{thm}
\label{pr-tf-equiv}
$\P^-_r\Lambda^k(T,f) = E[\P^-_r\Lambda^k(f)]$ and 
for $\mu\in\P^-_r\Lambda^k(f)$, 
$E \mu= E_{f,T}^{k,r,-}\mu$ can be characterized as the unique extension
of $\mu$ to $\P^-_r\Lambda^k(T,f)$.
\end{thm}

%To prove this theorem, we first establish the following result.
%\begin{lem}\label{P--space-prop}
%Let $f \in \Delta(T)$ and $f^* = \{x_0\}$, where $x_0$ is a vertex of $T$.
%Then $\P^-_r\Lambda^k(T,f)= E[\P^-_r\Lambda^k(f)]$.
%\end{lem}

\begin{proof}
We note that the second statement of the theorem follows from the first, since
$\tr_{T,f}$ from $E[P^-_r\Lambda^k(f)]$ to $P^-_r\Lambda^k(f)$ has a unique
right inverse.  Since $E[P^-_r\Lambda^k(f)] \subseteq
P^-_r\Lambda^k(T,f)$, we only need to prove the opposite inclusion.  Without
loss of generality we may assume that $f=[x_{m+1},\ldots,x_n]$, $f^* =
[x_0,\ldots,x_m]$, for some $0\le m< n$.  We proceed by induction on $m$.
When $m=0$, we may assume without loss of generality 
that the vertex $x_0$ is at the
origin.  Now $\P^-_r\Lambda^k=\P_{r-1}\Lambda^k + \kappa
\H_{r-1}\Lambda^{k+1}$ (where $\H_{r}$ denotes the homogeneous polynomials of
the degree $r$).  Since $\kappa \H_{r-1}\Lambda^{k+1}\subseteq\H_r\Lambda^k$,
every element $\omega\in \kappa \H_{r-1}\Lambda^{k+1}$ vanishes to order $r$
at the origin.  On the other hand, no non-zero element of $\P_{r-1}\Lambda^k$
vanishes to order $r$ at the origin.  Thus $\P^-_r\Lambda^k(T,f) = \kappa
\H_{r-1}\Lambda^{k+1}$.  It
follows from \cite[Theorem 3.3]{acta} that
\[
\dim \P^-_r\Lambda^k(T,f) =
\dim \kappa \H_{r-1}\Lambda^{k+1}(T) = 
\binom{r+n-1}{n-k-1}\binom{r+k-1}{k} = \dim \P^-_r\Lambda^k(f)
= \dim E[P^-_r\Lambda^k(f)],
\]
and since $E[P^-_r\Lambda^k(f)] \subseteq \P^-_r\Lambda^k(T,f)$,
$E[P^-_r\Lambda^k(f)] = \P^-_r\Lambda^k(T,f)$.

Now suppose that $\omega$ vanishes to order $r$ on the $m$-dimensional face
$[x_0,\ldots,x_m]$ with $m>0$.  Let $T'=[x_1,\ldots,x_n]$,
$\omega'=\tr_{T,T'}\omega$.  Then $\omega'\in\P^-_r\Lambda^k(T')$ vanishes to
order $r$ on the $(m-1)$-dimensional face $[x_1,\ldots,x_m]$, so, by
induction, $\omega' = E_{f,T'} \mu$  for some $\mu \in P^-_r\Lambda^k(f)$.
Furthermore, since $\omega$ vanishes to order $r$ at $x_0$, we
can use the result established above for $m=0$ to conclude that $\omega =
E_{T',T} \omega' = E_{T',T} E_{f,T'} \mu$.  However, it follows immediately
from \eqref{bary-ext-} that $E_{T',T} E_{f,T'}  = E_{f,T}$,
and hence the two spaces are equal, and the theorem is established.
\end{proof}

\section{A geometric decomposition of $\P_r\Lambda^k(T)$}\label{sec:decomp}

In this section, we again apply the theory developed in Section~\ref{cgd},
this time with $X(T) = \P_r\Lambda^k(T)$.  
In this case condition \eqref{cassump} is obvious, since
$\tr_{T,f}\P_r\Lambda^k(T)= \P_r\Lambda^k(f)$.
In view
of the previous section, one might hope that we could define the extension operator as
\begin{equation*}
(\lambda^f)^{\alpha} d\lambda^f_{\sigma} \mapsto
\lambda^{\alpha} d \lambda _{\sigma}.
\end{equation*}
However, this does not lead to a well-defined operator.  To appreciate the problem,
consider the space $\P_2\Lambda^1(T)$, where $T
\subset \R^2$ is a triangle spanned by the vertices $x_0,x_1,x_2$, and let $f
= [x_1,x_2]$. Then $\lambda^f_1 \lambda^f_2 (d \lambda^f_1 + d
\lambda^f_2) =0$, but $\lambda_1 \lambda_2 (d
\lambda_1 + d \lambda_2) = - \lambda_1 \lambda_2 d \lambda_0 \neq 0$.

%Assumption~\ref{A0} would then follow easily using the spanning set
%\eqref{pr-span}, since the trace of a spanning function 
%$\lambda^\alpha d\lambda_\sigma$ is the analogous spanning function
%$(\lambda^f)^\alpha d\lambda^f_\sigma$ for $\P_r \Lambda^k(f)$ if
%$\llbracket\alpha,\sigma\rrbracket \subset \I(f)$, and Assumption~\ref{A1}
%follows immediately from the definition of $\P_r \Lambda^k(T,f)$.
%The problem with this definition is that Assumption~\ref{A2} does not hold.
%For example,  consider the space $\P_2\Lambda^1(S)$, where $S
%\subset \R^2$ is a triangle spanned by the vertices $x_0,x_1,x_2$, and let $f
%= [x_1,x_2]$. Then $\tr_{S,f} \lambda_1 \lambda_2 (d
%\lambda_1 + d \lambda_2) = \lambda^f_1 \lambda^f_2 (d \lambda^f_1 + d
%\lambda^f_2) =0$, but $\lambda_1 \lambda_2 (d
%\lambda_1 + d \lambda_2) = - \lambda_1 \lambda_2 d \lambda_0 \neq 0$.

To remedy this situation, we will show that for $f,g \in \Delta(T)$ with $f \subseteq g$, a
consistent extension operator $E =E_{f,g}^{k,r}: \P_r\Lambda^k(f) \to
\P_r\Lambda^k(g)$ is given by
\begin{equation}
\label{bary-ext}
(\lambda^f)^{\alpha} d\lambda^f_{\sigma} \mapsto
(\lambda^g)^{\alpha} \psi^{\alpha,f,g}_{\sigma}, \qquad
\llbracket\alpha,\sigma\rrbracket \subseteq \I(f),
\end{equation}
where $\psi^{\alpha,f,g}_{\sigma}$ is defined as follows.
We first introduce forms $\psi_i^{\alpha,f,g} \in \Alt^1 T_g$ defined by
\begin{equation}
\label{psii-def}\psi_i^{\alpha,f,g} = d \lambda^g_i 
- \frac{\alpha_i}{|\alpha|} \sum_{j \in \I(f)} d \lambda^g_j,
 \quad i \in \I(f),
\end{equation}
and then define $\psi^{\alpha,f,g}_{\sigma} \in \Alt^k T_g$ by
\begin{equation}
\label{psisigma-def}
\psi^{\alpha,f,g}_{\sigma} = \psi^{\alpha,f,g}_{\sigma(1)} \wedge \cdots
\wedge \psi^{\alpha,f,g}_{\sigma(k)}, \quad \sigma\in\Sigma(1:k,0:n),\quad \llbracket\sigma\rrbracket\subseteq\I(f).
\end{equation}
A geometric interpretation of $\psi^{\alpha,f,g}_{\sigma}$ will be
given below.

First we show that $E$ is well-defined and is, in fact, an
extension operator.
\begin{thm}\label{P-spaces-tr}
There is a unique mapping $E=E_{f,g}^{k,r}: \P_r\Lambda^k(f) \to
\P_r\Lambda^k(g)$ satisfying \eqref{bary-ext}.  Moreover it is an
extension operator: $\tr_{g,f}E_{f,g}^{k,r}\omega=\omega$ for $\omega\in \P_r\Lambda^k(f)$.
\end{thm}
\begin{proof}
By the first part of Theorem~\ref{summary}, the set
\begin{equation*}
\{\,(\lambda^f)^\alpha d\lambda^f_\sigma\,|\,\alpha\in\N_0^{0:n},|\alpha|=r,
\sigma\in\Sigma(1:k,0:n),\ \llbracket\alpha, \sigma\rrbracket\subseteq \I(f),
\ \min\llbracket\sigma\rrbracket>\min\I(f) \,\}
\end{equation*}
is a basis for $\P_r\Lambda^k(f)$. Hence, we can define an extension
$E$ by \eqref{bary-ext}, if we restrict to the basis functions. We now show that
\eqref{bary-ext} holds for all $(\lambda^f)^\alpha d\lambda^f_\sigma$ with
$\llbracket\alpha, \sigma\rrbracket \subseteq \I(f)$, i.e., also when
$\min\llbracket\sigma\rrbracket=\min\I(f)$.
Writing $d \lambda_\sigma = d \lambda_{\sigma(1)} \wedge d \lambda_{\rho}$,
and using the fact that $\sum_{j \in \I(f)} d \lambda^f_j =0$ on the 
face $f$, we can write
\begin{equation*}
d \lambda^f_\sigma = - \sum_{\substack{j \in \I(f) \\ j \neq \sigma(1)}}
d \lambda^f_j \wedge d \lambda^f_{\rho}.
\end{equation*}
Hence,
\begin{equation*}
E[(\lambda^f)^{\alpha} d\lambda^f_{\sigma}] 
= - E[ (\lambda^f)^{\alpha} \sum_{\substack{j \in \I(f) \\ j \neq \sigma(1)}}
d \lambda^f_j \wedge d \lambda^f_{\rho}]
= - (\lambda^g)^{\alpha} \sum_{\substack{j \in \I(f) \\ j \neq \sigma(1)}}
\psi^{\alpha, f,g}_j \wedge \psi^{\alpha, f,g}_{\rho}
= (\lambda^g)^\alpha \psi^{\alpha, f,g}_\sigma,
\end{equation*}
where  in the last step we have used the fact that $\sum_{i \in \I(f)} \psi^{\alpha,f,g}_i =0$.

That $E$ is an extension operator follows directly from the observation
\begin{equation*}
\tr_{T,f} \psi_{\sigma} = d\lambda^f_{\sigma},
\end{equation*}
which holds since $\sum_{j \in \I(f)} d \lambda^f_j =0$ on the face $f$.
\end{proof}

\begin{thm}
\label{Econsistent}
The family of extension operators $E$ is consistent, i.e., for all $f,g,h\in\Delta(\T)$
with $f,g \subseteq h$ and all
$\omega \in \P_r\Lambda^k(f)$,
\begin{equation*}
\tr_{h,g}E_{f,h} \omega = E_{f \cap g, g} \tr_{f,f \cap g}\omega.
\end{equation*}
\end{thm}
\begin{proof}
It is enough to establish this result for $\omega = (\lambda^f)^{\alpha} 
d\lambda^f_{\sigma}$, with $\llbracket\alpha,\sigma\rrbracket\subseteq\I(f)$.
Now for such pairs $(\alpha, \sigma)$, $E_{f,h}[(\lambda^f)^{\alpha} (d
\lambda^f)_{\sigma}] = (\lambda^h)^{\alpha} \psi^{\alpha,f,h}_{\sigma}$.  To
determine $\tr_{h,g}[(\lambda^h)^{\alpha} \psi^{\alpha,f,h}_{\sigma}]$, we
consider three cases.  When $\llbracket\alpha\rrbracket\subseteq \I(f)$, but
$\llbracket\alpha\rrbracket\nsubseteq \I(g)$, $\tr_{h,g}[(\lambda^h)^{\alpha}
\psi^{\alpha,f,h}_{\sigma}] =0$, since
$\tr_{h,g}[(\lambda^h)^{\alpha}]=0$. If $\llbracket\alpha\rrbracket\subseteq
\I(f \cap g)$, then $\tr_{h,g}[(\lambda^h)^{\alpha}] = (\lambda^g)^{\alpha}$,
so we need only compute $\tr_{h,g} \psi^{\alpha,f,h}_{\sigma}$. We do
this by first considering $\tr_{h,g} \psi^{\alpha,f,h}_{i}$ for
$i \in \I(f)$. If $i \in \I(f \cap g)$, we have
\begin{equation*}
\tr_{h,g} \psi^{\alpha,f,h}_{i} = \tr_{h,g}(
d \lambda^h_i - \frac{\alpha_i}{|\alpha|} \sum_{j \in \I(f)}
d \lambda^h_j)
= d \lambda^g_i - \frac{\alpha_i}{|\alpha|} \sum_{j \in \I(f \cap g)}
d \lambda^g_j = \psi^{\alpha,f \cap g,g}_{i}.
\end{equation*}
On the other hand, if $i \in \I(f) \setminus \I(f \cap g)$, then
since $\alpha_i =0$, $\psi^{\alpha,f,h}_{i} = d \lambda^h_i$ and so
$\tr_{h,g} \psi^{\alpha,f,h}_{i} = 0$.
Combining these results, we obtain
\begin{equation*}
\tr_{h,g}E_{f,h}[(\lambda^f)^{\alpha} d\lambda^f_{\sigma}] = \begin{cases}
(\lambda^g)^{\alpha} \psi^{\alpha, f \cap g,g}_{\sigma}, & \text{if } 
\llbracket\alpha,\sigma\rrbracket\subseteq\I(f \cap g),
\\ 0, & \text{otherwise}. \end{cases}
\end{equation*}
But
\begin{equation*}
\tr_{f,f \cap g} [(\lambda^f)^{\alpha} d\lambda^f_{\sigma}]
= \begin{cases}
(\lambda^{f \cap g})^{\alpha} d \lambda^{f \cap g}_{\sigma}, & \text{if }
\llbracket\alpha,\sigma\rrbracket\subseteq\I(f \cap g), \\
0, & \text{otherwise}, \end{cases}
\end{equation*}
and hence
\begin{equation*}
E_{f \cap g, g} \tr_{f,f \cap g} [(\lambda^f)^{\alpha} d\lambda^f_{\sigma}]
= \begin{cases}
(\lambda^{g})^{\alpha} \psi^{\alpha, f \cap g,g}_{\sigma}, & \text{if }
\llbracket\alpha,\sigma\rrbracket\subseteq\I(f \cap g), \\
0, & \text{otherwise}. \end{cases}
\end{equation*}
\end{proof}

%\begin{remark}
%The forms $\psi_i$ that we have presented may still appear unmotivated.
%However, since we want property \eqref{note-sum}, the most natural
%modification of $d \lambda_i$ is to subtract a multiple of an expression that
%vanishes on the face $f$, which is what we have done.
%Note that the specific choice of
%the factor $\alpha_i/|\alpha|$ is used above in an essential way to show that
%the extension operator is consistent.
%$\Box$
%\end{remark}

From Theorem~\ref{basic-decomp}, we obtain the desired geometric decomposition
$\P_r\Lambda^k(\T)$.
\begin{thm}
\label{direct-sum-decomp}
\begin{equation*}
\P_r\Lambda^k(\T) = \bigoplus_{\substack{f\in\Delta(\T)
\\ \dim f\ge k}} E_f[\0\P_r\Lambda^k(f)].
\end{equation*}
where $E_f:\0\P_r\Lambda^k(f)\to \P_r\Lambda^k(\T)$ denotes 
the global extension
operator determined by the family $E_{f,g}^{k,r}$.
\end{thm}

Combining this result with the second part of Theorem~\ref{summary}, we obtain
an explicit spanning set and basis for $\P^-_r\Lambda^k(\T)$ (see 
Section~\ref{sec:bases}).

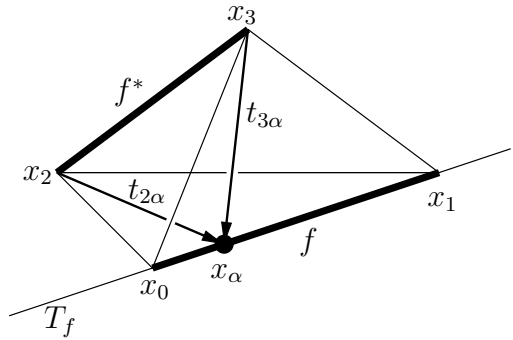
\begin{figure}[ht]
\begin{picture}(0,0)%
\includegraphics{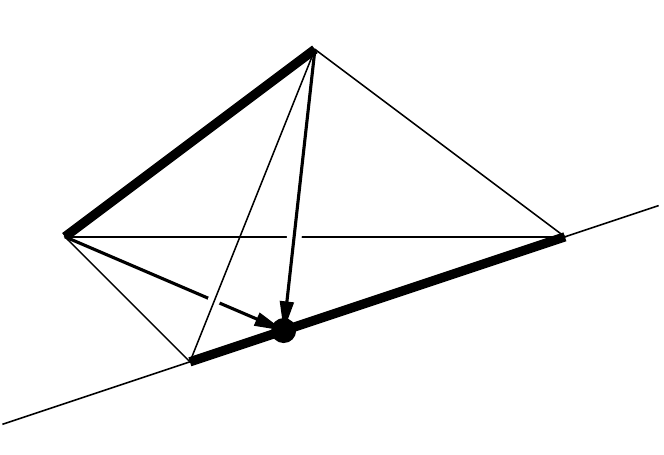}%
\end{picture}%
\setlength{\unitlength}{3947sp}%
\begingroup\makeatletter\ifx\SetFigFontNFSS\undefined%
\gdef\SetFigFontNFSS#1#2#3#4#5{%
  \reset@font\fontsize{#1}{#2pt}%
  \fontfamily{#3}\fontseries{#4}\fontshape{#5}%
  \selectfont}%
\fi\endgroup%
\begin{picture}(3174,2167)(2089,-1107)
\put(2321,-1038){\makebox(0,0)[lb]{\smash{{\SetFigFontNFSS{12}{14.4}{\familydefault}{\mddefault}{\updefault}{\color[rgb]{0,0,0}$T_f$}%
}}}}
\put(3921,-544){\makebox(0,0)[lb]{\smash{{\SetFigFontNFSS{12}{14.4}{\familydefault}{\mddefault}{\updefault}{\color[rgb]{0,0,0}$f$}%
}}}}
\put(2921,-823){\makebox(0,0)[lb]{\smash{{\SetFigFontNFSS{12}{14.4}{\familydefault}{\mddefault}{\updefault}{\color[rgb]{0,0,0}$x_0$}%
}}}}
\put(2848,-224){\makebox(0,0)[lb]{\smash{{\SetFigFontNFSS{12}{14.4}{\familydefault}{\mddefault}{\updefault}{\color[rgb]{0,0,0}$t_{2\alpha}$}%
}}}}
\put(4728,-270){\makebox(0,0)[lb]{\smash{{\SetFigFontNFSS{12}{14.4}{\familydefault}{\mddefault}{\updefault}{\color[rgb]{0,0,0}$x_1$}%
}}}}
\put(3470,901){\makebox(0,0)[lb]{\smash{{\SetFigFontNFSS{12}{14.4}{\familydefault}{\mddefault}{\updefault}{\color[rgb]{0,0,0}$x_3$}%
}}}}
\put(2187,-104){\makebox(0,0)[lb]{\smash{{\SetFigFontNFSS{12}{14.4}{\familydefault}{\mddefault}{\updefault}{\color[rgb]{0,0,0}$x_2$}%
}}}}
\put(2757,395){\makebox(0,0)[lb]{\smash{{\SetFigFontNFSS{12}{14.4}{\familydefault}{\mddefault}{\updefault}{\color[rgb]{0,0,0}$f^*$}%
}}}}
\put(3366,-720){\makebox(0,0)[lb]{\smash{{\SetFigFontNFSS{12}{14.4}{\familydefault}{\mddefault}{\updefault}{\color[rgb]{0,0,0}$x_\alpha$}%
}}}}
\put(3589,259){\makebox(0,0)[lb]{\smash{{\SetFigFontNFSS{12}{14.4}{\familydefault}{\mddefault}{\updefault}{\color[rgb]{0,0,0}$t_{3\alpha}$}%
}}}}
\end{picture}%
\caption{$T=[x_0,x_1,x_2,x_3]$, $f=[x_0,x_1]$, $\alpha=(3,1,0,0)$, $\R^3= T_f \oplus \spn\{t_{2\alpha},t_{3\alpha}\}$.}
\label{fg:projdiag}
\end{figure}

We now turn to a geometric characterization of the extension operator
$E = E_{f,g}^{k,r}:\P_r\Lambda^k(f)\to\P_r\Lambda^k(T,f)$.  
First, we will motivate the choice of $E$, and in particular
the forms $\psi_{\sigma}^{\alpha,f,g}$,  by establishing some 
some additional properties of these forms.
Observe that any multi-index $\alpha$ determines a convex combination of the
vertices $x_i$ of $T$, namely
\[
x_\alpha:=|\alpha|^{-1}\sum_m\alpha_m x_m\in T, 
\]
and if $\llbracket\alpha\rrbracket\subset\I(f)$, then $x_\alpha\in f$.
For each such multi-index $\alpha$, we then define the vectors
\begin{equation*}
t_{\alpha l} = x_\alpha - x_l = \frac1{|\alpha|}\sum_{m \in \I(f)} \alpha_m (x_m-x_l), \quad l \in \I(f^*).
\end{equation*}
Clearly, for each such $\alpha$,
$\R^n$ decomposes as  the direct sum
$T_f \oplus \spn\{t_{\alpha l} \, | \, l \in \I(f^*)\}$,
where $T_f$ denotes the tangent space of $f$.  See Figure~\ref{fg:projdiag}.
This decomposition
defines a projection operator $P=\Pfalpha : \R^n \to T_f$
determined by the equations
$P v =v$ for $v \in T_f$ and $P \, t_{\alpha l} = 0$ for $l \in \I(f^*)$.
Hence, we have
\begin{equation}
\label{chara}
P_{f,\alpha}^* \Alt^k T_f = \{a \in  \Alt^k \R^n \, | \, a \lrcorner
\, t_{\alpha l} =0, \, l \in \I(f^*)\}.
\end{equation}
Furthermore, since $d \lambda_j(x_m-x_l) = \delta_{jm}$
for any $j \in \I(f)$, $m \in \I(f)$, and $l \in \I(f^*)$,
we get for 
$\llbracket\alpha\rrbracket \subset \I(f)$,
\begin{equation}
\label{psii-contract}
\psi_i^{\alpha,f,T}(t_{\alpha l}) = d \lambda_i (t_{\alpha l})
- \frac{\alpha_i}{|\alpha|} \sum_{j \in \I(f)} d \lambda_j(t_{\alpha l})
= \frac1{\alpha}\left(\alpha_i - \frac{\alpha_i}{|\alpha|} \sum_{j \in \I(f)} \alpha_j \right)=0.
\end{equation}
It follows that
\begin{equation}\label{psitrace}
\psi_{\sigma}^{\alpha,f,T}(v_1, \cdots, v_k) = 
\psi_{\sigma}^{\alpha,f,T}(P v_1, \cdots, P v_k) 
= \tr_{T,f}\psi_{\sigma}^{\alpha,f,T}(P v_1, \cdots, P v_k) 
= d\lambda^f_{\sigma} (P v_1, \cdots, P v_k).
\end{equation}
Hence, in the language of pullbacks,
\begin{equation*}
\psi_{\sigma}^{\alpha,f,T} =
P_{f, \alpha}^* d\lambda^f_{\sigma},
\end{equation*}
where $P_{f, \alpha}^*$ is the pullback of $P_{f, \alpha}$, and so
\begin{equation*}
E_{f,g}[(\lambda^f)^{\alpha} d\lambda^f_{\sigma}]
= (\lambda^g)^{\alpha} P_{f,\alpha}^* d\lambda^f_{\sigma}.
\end{equation*}

Recall that the geometric characterization in the previous section hinged upon
the fact that a form in $\P_r^- \Lambda^k(T)$ which vanishes to order $r$ on
$f^*$ and has vanishing trace on $f$ must vanish identically.  Now this is not
true for an arbitrary element of the larger space $\P_r \Lambda^k(T)$.
Returning to the example given at the beginning of this section, where
$T=[x_0,x_1,x_2]$ and $f = [x_1,x_2]$, the form
\begin{equation*}
\omega = \lambda_1 \lambda_2 [d \lambda_1 + d \lambda_2] 
= - \lambda_1 \lambda_2 d \lambda_0 \in \P_2 \Lambda^1(T)
\end{equation*}
vanishes to second order at $f^* = \{x_0 \}$. However, $\tr_{T,f} \omega$ also
vanishes. Thus, additional conditions on $\omega$ will be needed in order to
insure that $\omega$ is uniquely determined by $\tr_{T,f} \omega$.
We say that {\em $\omega$ vanishes to order $r^+$ on $f^*$},
if $\omega$ vanishes to order $r$  on $f^*$ and the following conditions hold:
\begin{equation}
\label{gen-x-cond}
\partial_{t_l}^\alpha \omega \lrcorner t_{\alpha l} = 0, \  l \in
\I(f^*), \  \llbracket\alpha \rrbracket \subseteq \I(f), \  |\alpha| = r.
\end{equation}
Here $\partial^\alpha_{t_l}
:= \prod_{j \in \I(f)} \partial_{t_{jl}}^{\alpha_j}$ with $\partial_{t_{jl}}=
t_{jl}\cdot\nabla$,  the directional derivative along the vector
$t_{jl}:=x_j-x_l$.  The contraction operator $\lrcorner$ is defined
at the start of Section~\ref{df}.

Note that $\partial_{t_{jl}}
\lambda_i = \delta_{ij}$ for $i,j\in\I(f)$, $j\in\I(f^*)$.   It follows that if $\alpha, \beta \in \N_0^n$,
with $|\alpha| = |\beta| =r$, $\llbracket\alpha\rrbracket,\llbracket\beta\rrbracket\subseteq\I(f)$, and $l\in\I(f^*)$, then
\begin{equation}
\label{delta-prop}
\partial_{t_l}^{\beta} \lambda^{\alpha} = 0 \quad \text{for } \alpha \neq \beta
 \text{ and } \partial_{t_l}^{\alpha} \lambda^{\alpha} = \alpha!.
\end{equation}

Setting
\begin{equation*}
\P_r\Lambda^k(T,f) = \{\, \omega \in \P_r\Lambda^k(T) \,|\,
\omega \text{ vanishes to order } r^+ \text{ on $f^*$} \, \},
\end{equation*}
we can now give the geometric description of the extension operator $E$.
\begin{thm}
\label{prtf-equiv}
$\P_r\Lambda^k(T,f) = E[\P_r\Lambda^k(f)]$ and 
for $\mu\in\P_r\Lambda^k(f)$, 
$E \mu= E_{f,T}^{k,r}\mu$ can be characterized as the unique extension
of $\mu$ to $\P_r\Lambda^k(T,f)$.
\end{thm}
\begin{proof}
We note that the second statement of the theorem follows from the first, since
$\tr_{T,f}$ from $E[P_r\Lambda^k(f)]$ to $P_r\Lambda^k(f)$ has a unique right
inverse. To prove the first statement, we first show that $E[\P_r\Lambda^k(f)]
\subseteq \P_r\Lambda^k(T,f)$.  Observe first that $E[(\lambda^f)^{\alpha} (d
\lambda^f)_{\sigma}] = \lambda^{\alpha} \psi^{\alpha,f,T}_\sigma$ vanishes
to order $r$ on $f^*$ since $(\lambda)^{\alpha}$ does.  Next, note that
\eqref{delta-prop} tells us that
$\partial^\beta_{t_l}[\lambda^\alpha\psi^{\alpha,f,T}_\sigma] = \alpha!
\psi^{\alpha,f,T}_\sigma$ if $\beta=\alpha$ and vanishes if $\beta$ is any
other multi-index of order $r$ with $\llbracket\beta \rrbracket \subseteq
\I(f)$.  Therefore, the conditions \eqref{gen-x-cond} for vanishing of order
$r^+$ are reduced to verifying the conditions $\psi^{\alpha,f,T}_\sigma
\lrcorner t_{\alpha l} = 0$ for all $l \in \I(f^*)$. However, this follows
immediately from the definition of the wedge product and
\eqref{psii-contract}.

To show that $\P_r\Lambda^k(T,f) \subseteq E[\P_r\Lambda^k(f)]$, we
use Lemma~\ref{0-forms} to see that any element $\omega \in
\P_r\Lambda^k(T,f)$ admits a representation of the form
\begin{equation*}
\omega = \sum_{\substack{\llbracket \alpha \rrbracket \subseteq
    \I(f) \\
    |\alpha|=r}} a_\alpha
\lambda^\alpha
\end{equation*}
for some $a_\alpha \in \Alt^k \R^n$.  However,
invoking \eqref{delta-prop} and \eqref{gen-x-cond}, we conclude that, if
$\omega$ vanishes to the order $r^+$ on $f^*$, then $a_\alpha\lrcorner
t_{\alpha l} = 0$ for all $l \in \I(f^*)$, and hence by \eqref{chara},
$a_{\alpha} \in P_{f,\alpha}^* \Alt^k T_f$.  It therefore follows
from \eqref{psitrace}, that $\omega  \in E[\P_r\Lambda^k(f)]$.
\end{proof}

\section{Construction of bases}
\label{sec:bases}

From Theorem~\ref{direct-sum-decomp-}, \eqref{eft},
Theorem~\ref{P--spaces-tr}, and part 4 of Theorem~\ref{summary}, one
immediately obtains explicit formulas for a spanning set and basis for
$\P^-_r\Lambda^k(\T)$, with each spanning and basis form associated to a particular face $f
\in \Delta(\T)$.  The forms associated to $f$ vanish on simplices $T \in \T$ that do
not contain $f$, while for $T$ containing $f$, the spanning and basis forms are given by
\begin{equation*}
\{\,(\lambda^T)^\alpha\phi^T_\sigma\,|\,
\alpha\in\mathbb N_0^{0:n}, |\alpha|=r-1,\sigma\in\Sigma(0:k,0:n),
\llbracket\alpha,\sigma\rrbracket = \I(f) \,\}
\end{equation*}
and
\begin{equation*}
\{\,(\lambda^T)^\alpha\phi^T_\sigma\,|\,
\alpha\in\mathbb N_0^{0:n}, |\alpha|=r-1,\sigma\in\Sigma(0:k,0:n),
\llbracket\alpha,\sigma\rrbracket = \I(f), \alpha_i =0
\text{ if } i < \min \llbracket\sigma\rrbracket \,\},
\end{equation*}
respectively.  Note that the spanning set is independent of the ordering of
the vertices, while our choice of basis depends on the ordering of the
vertices.  Other choices of basis are possible as well, but there is no one
canonical choice.

The same considerations give an explicit spanning set and basis for
$\P_r\Lambda^k(\T)$, based on Theorem~\ref{direct-sum-decomp},
\eqref{eft}, Theorem~\ref{P-spaces-tr}, and
part 2 of Theorem~\ref{summary}.  The corresponding formulas for the spanning
set and basis are:
\begin{equation*}
\{\,(\lambda^T)^\alpha \psi^{\alpha,f,T}_\sigma\,|\,
\alpha\in\mathbb N_0^{0:n}, |\alpha|=r,\sigma\in\Sigma(1:k,0:n),
%\llbracket\alpha,\sigma\rrbracket = \I(f) \,\}.
\end{equation*}
and
\begin{equation*}
\{\,(\lambda^T)^\alpha \psi^{\alpha,f,T}_\sigma\,|\,
\alpha\in\mathbb N_0^{0:n}, |\alpha|=r,\sigma\in\Sigma(1:k,0:n),
\llbracket\alpha,\sigma\rrbracket = \I(f),\ \alpha_i=0\text{ if }
i<\min[I(f)\setminus \llbracket\sigma\rrbracket \,]\,\},
\end{equation*}
respectively, where $\psi^{\alpha,f,T}_{\sigma}$
is defined by \eqref{psii-def} and \eqref{psisigma-def}.

%\begin{equation*}
%\psi_i = \psi_i^{\alpha,f,T} = d \lambda^T_i
%- \frac{\alpha_i}{|\alpha|} \sum_{j \in \I(f)} d \lambda^T_j,
%\quad i \in \I(f),
%\end{equation*}
%and
%\begin{equation*}
%\psi^T_{\sigma} = \psi^{\alpha,f,T}_{\sigma}
%= \psi_{\sigma(1)} \wedge \cdots
%\wedge \psi_{\sigma(k)}, \quad \sigma\in\Sigma(1:k,0:n),
%\quad \llbracket\sigma\rrbracket\subseteq\I(f).
%\end{equation*}

Bases for the spaces $\P^-_r \Lambda^k$ and $\P_r \Lambda^k$ are summarized in
Tables~\ref{tb:t1}--\ref{tb:t4} for $n=2,3$, $0 < k < n$, $r=1,2,3$. In the tables, we
assume $i < j < k < l$, and recall that the Whitney forms $\phi_{ij}$ and
$\phi_{ijk}$ are given by:
\begin{equation*}
\phi_{ij} = \lambda_i \,d \lambda_j - \lambda_j\, d \lambda_i, \qquad
\phi_{ijk} = \lambda_i \,d \lambda_j \wedge d \lambda_k
- \lambda_j \,d \lambda_i \wedge d \lambda_k
+ \lambda_k \,d \lambda_i \wedge d \lambda_j.
\end{equation*}

\begin{table}[htb]
\caption{Bases for the spaces $\P^-_r \Lambda^1$ and
$\P_r \Lambda^1$, $n =2$.}
\label{tb:t1}
\footnotesize
\begin{center}
\begin{tabular}{ccccccc}
\hline
\\
\multicolumn{1}{c}{} && \multicolumn{2}{c}{$\P^-_r \Lambda^1$}
&& \multicolumn{2}{c}{$\P_r \Lambda^1$} \\
 && \multicolumn{2}{c}{\rule[5pt]{2in}{.5pt}}
&& \multicolumn{2}{c}{\rule[5pt]{2in}{.5pt}}
%\hline
\\
$r$ && Edge $[x_i,x_j]$ & Triangle $[x_i,x_j,x_k]$ &&  Edge $[x_i,x_j]$
& Triangle $[x_i,x_j,x_k]$ \\
\hline
$1$\rule{0pt}{15pt} && $\phi_{ij}$  &   && $\lambda_i d \lambda_j$,  \
$\lambda_j d \lambda_i$ & \\[2ex]
%\hline
$2$ && $\{\lambda_i, \lambda_j\} \phi_{ij}$ &
$\lambda_k \phi_{ij}$, \ $\lambda_j \phi_{ik}$  &&
$\lambda_i^2 d \lambda_j$, \ $\lambda_j^2 d \lambda_i$
& $\lambda_i \lambda_j d \lambda_k$, \ $\lambda_i \lambda_k d \lambda_j$
\\
&&  &   && $\lambda_i \lambda_j d(\lambda_j - \lambda_i)$
& $\lambda_j \lambda_k d \lambda_i$
\\[2ex]
%\hline
$3$ && $\{\lambda_i^2, \lambda_j^2, \lambda_i \lambda_j\} \phi_{ij}$ &
$\{\lambda_i, \lambda_j, \lambda_k\} \lambda_k \phi_{ij}$
&&
$\lambda_i^3 d \lambda_j$, \ $\lambda_j^3 d \lambda_i$
&
$\{\lambda_i, \lambda_j, \lambda_k\} \lambda_i \lambda_j d \lambda_k$,

\\
&&  &
$\{\lambda_i, \lambda_j, \lambda_k\} \lambda_j \phi_{ik}$
&& $\lambda_i^2 \lambda_j d(2 \lambda_j - \lambda_i)$ &
$\{\lambda_i, \lambda_j, \lambda_k\} \lambda_i \lambda_k d \lambda_j$
\\
&& &  && $\lambda_i \lambda_j^2 d(\lambda_j - 2 \lambda_i)$ &
$\{\lambda_j, \lambda_k\} \lambda_j \lambda_k d \lambda_i$
\\
\\
\hline
\end{tabular}
\end{center}
\end{table}

\begin{table}[htb]
\caption{Bases for the spaces $\P^-_r \Lambda^1$ and
$\P^-_r \Lambda^2$, $n =3$.}
\label{tb:t2}
\footnotesize
\begin{center}
\begin{tabular}{ccccccc}
\hline
\\
\multicolumn{1}{c}{} & \multicolumn{3}{c}{$\P^-_r \Lambda^1$}
&& \multicolumn{2}{c}{$\P^-_r \Lambda^2$} \\
%\hline
 & \multicolumn{3}{c}{\rule[5pt]{2.9in}{.5pt}}
&& \multicolumn{2}{c}{\rule[5pt]{2.2in}{.5pt}}
\\
$r$ & Edge $[x_i,x_j]$ & Face $[x_i,x_j,x_k]$ & Tet $[x_i,x_j,x_k,x_l]$
&& Face $[x_i,x_j,x_k]$ & Tet $[x_i,x_j,x_k,x_l]$ \\
\hline
$1$\rule{0pt}{15pt} & $\phi_{ij}$ &   & & & $\phi_{ijk}$ & \\[2ex]
%\hline
$2$ & $\{\lambda_i, \lambda_j\} \phi_{ij}$  &
$\lambda_k \phi_{ij}$,  \ $\lambda_j \phi_{ik}$ & & &
$\{\lambda_i,\lambda_j, \lambda_k\} \phi_{ijk}$ &
$\lambda_l \phi_{ijk}$, \ $\lambda_k \phi_{ijl}$
\\
&   &
&  &  &  & $\lambda_j \phi_{ikl}$
\\[2ex]
%\hline
$3$ & $\{\lambda_i^2, \lambda_j^2, \lambda_i \lambda_j\} \phi_{ij}$ &
$\{\lambda_i, \lambda_j, \lambda_k\} \lambda_k \phi_{ij}$
& $\lambda_k \lambda_l \phi_{ij}$ &&
$\{\lambda_i^2, \lambda_j^2, \lambda_k^2\} \phi_{ijk}$
& $\{\lambda_i, \lambda_j, \lambda_k, \lambda_l\} \lambda_l \phi_{ijk}$
\\
&   & $\{\lambda_i, \lambda_j, \lambda_k\} \lambda_j \phi_{ik}$
& $\lambda_j \lambda_l \phi_{ik}$
& & $\{\lambda_i \lambda_j, \lambda_i \lambda_k,
\lambda_j \lambda_k\} \phi_{ijk}$ &
$\{\lambda_i, \lambda_j, \lambda_k, \lambda_l\} \lambda_k \phi_{ijl}$
\\
&  &    & $\lambda_j \lambda_k \phi_{il}$
 && &
$\{\lambda_i, \lambda_j, \lambda_k, \lambda_l\} \lambda_j \phi_{ikl}$
\\
\\
\hline
\end{tabular}
\end{center}
\end{table}

\begin{table}[htb]
\caption{Basis for the space $\P_r \Lambda^1$, $n = 3$.}
\label{tb:t3}
\footnotesize
\begin{center}
\begin{tabular}{c c c c c c c}
%\hline
%\multicolumn{1}{c}{} & \multicolumn{3}{c}{$\P_r \Lambda^1$} \\
\hline
\\
$r$ &\hspace{.25in}& Edge $[x_i,x_j]$ &\hspace{.25in}& Face $[x_i,x_j,x_k]$ 
&\hspace{.25in}& Tet $[x_i,x_j,x_k,x_l]$
\\
\hline
$1$\rule{0pt}{15pt} && $\lambda_i d \lambda_j$, \ $\lambda_j d \lambda_i$
&& && \\[2ex]
%\hline
$2$ &&
$\lambda_i^2 d \lambda_j$, \ $\lambda_j^2 d \lambda_i$, \
$\lambda_i \lambda_j d(\lambda_j - \lambda_i)$
&& $\lambda_i \lambda_j d \lambda_k$, \ $\lambda_i \lambda_k d \lambda_j$, \
$\lambda_j \lambda_k d \lambda_i$ &&
\\[2ex]
%\hline
$3$ && $\lambda_i^3 d \lambda_j$, \
$\lambda_i^2 \lambda_j d(2 \lambda_j - \lambda_i)$
&&
$\{\lambda_i, \lambda_j\} \lambda_i \lambda_j d \lambda_k$, \
$\lambda_i \lambda_j \lambda_k d (2\lambda_k - \lambda_i - \lambda_j)$
&& $\lambda_i \lambda_j \lambda_k d \lambda_l$, \
$\lambda_i \lambda_j \lambda_l d \lambda_k$
\\
&& $\lambda_j^3 d \lambda_i$, \
$\lambda_i \lambda_j^2  d(\lambda_j - 2 \lambda_i)$   &&
$\{\lambda_i, \lambda_k\} \lambda_i \lambda_k d \lambda_j$, \
$\lambda_i \lambda_j \lambda_k d (2\lambda_j - \lambda_i - \lambda_k)$
&&
$\lambda_i \lambda_k \lambda_l d \lambda_j$, \
$\lambda_j \lambda_k \lambda_l d \lambda_i$
\\
&& && $\{\lambda_j, \lambda_k\} \lambda_j \lambda_k d \lambda_i$ &&
\\
\\
\hline
\end{tabular}
\end{center}
\end{table}

\begin{table}[htb]
\caption{Basis for the space $\P_r \Lambda^2$, $n = 3$.}
\label{tb:t4}
\footnotesize
\begin{center}
\begin{tabular}{c c c c c}
%\hline
%\multicolumn{1}{} & \multicolumn{2}{c}{$\P_r \Lambda^2$, $n = 3$.} \\
\hline
\\
$r$ &\hspace{.25in} & Face $[x_i,x_j,x_k]$ &\hspace{.25in} & Tet $[x_i,x_j,x_k,x_l]$ \\
\hline
$1$\rule{0pt}{15pt} && $\lambda_k d \lambda_i \wedge d \lambda_j$, \
$\lambda_j d \lambda_i \wedge d \lambda_k$, \
$\lambda_i d \lambda_j \wedge d \lambda_k$  &&
\\[2ex]
%\hline
$2$ && $\lambda_k^2 d \lambda_i \wedge d \lambda_j$, \
$\lambda_j \lambda_k d \lambda_i \wedge  d(\lambda_k - \lambda_j)$
&& $\lambda_k \lambda_l d \lambda_i \wedge d \lambda_j$, \
$\lambda_j \lambda_l d \lambda_i \wedge d \lambda_k$
\\
&&  $\lambda_j^2 d \lambda_i \wedge d \lambda_k$, \
$\lambda_i \lambda_j d(\lambda_j - \lambda_i) \wedge d \lambda_k$
 && $\lambda_j \lambda_k d \lambda_i \wedge d \lambda_l$, \
$\lambda_i \lambda_l d \lambda_j \wedge d \lambda_k$
\\
&&  $\lambda_i^2 d \lambda_j \wedge d \lambda_k$, \
 $\lambda_i \lambda_k d \lambda_j \wedge d (\lambda_k - \lambda_i)$
&& $\lambda_i \lambda_k d \lambda_j \wedge d \lambda_l$, \
$\lambda_i \lambda_j d \lambda_k \wedge d \lambda_l$
\\[2ex]
%\hline
$3$ && $\lambda_k^3 d \lambda_i \wedge d \lambda_j$, \
$\lambda_j^3 d \lambda_i \wedge d \lambda_k$, \
$\lambda_i^3 d \lambda_j \wedge d \lambda_k$
&& $\{\lambda_k, \lambda_l\} \lambda_k \lambda_l
d \lambda_i \wedge d \lambda_j$
\\
&&  $\lambda_j^2  \lambda_k d \lambda_i
\wedge d (2 \lambda_k - \lambda_j)$, \
$\lambda_j \lambda_k^2 d \lambda_i \wedge d (\lambda_k - 2 \lambda_j)$
&& $\{\lambda_j, \lambda_k, \lambda_l\} \lambda_j \lambda_l
d \lambda_i \wedge d \lambda_k$
\\
&& $\lambda_i^2 \lambda_j d(2 \lambda_j - \lambda_i) \wedge d \lambda_k$, \
$\lambda_i^2 \lambda_k d \lambda_j \wedge d(2\lambda_k - \lambda_i)$
&& $\{\lambda_j, \lambda_k, \lambda_l\} \lambda_j \lambda_k
d \lambda_i \wedge d \lambda_l$
\\
&& $\lambda_i \lambda_j^2 d(\lambda_j - 2 \lambda_i) \wedge d \lambda_k$, \
$\lambda_i \lambda_k^2 d \lambda_j \wedge d (\lambda_k - 2 \lambda_i)$
&&
$\{\lambda_i, \lambda_j, \lambda_k, \lambda_l\} \lambda_i \lambda_l
d \lambda_j \wedge d \lambda_k$
\\
&& $\lambda_i \lambda_j \lambda_k d (2 \lambda_j - \lambda_i - \lambda_k)
\wedge  d (2 \lambda_k - \lambda_i - \lambda_j)$ &&
$\{\lambda_i, \lambda_j, \lambda_k, \lambda_l\} \lambda_i \lambda_k
 d \lambda_j \wedge d \lambda_l$
\\
&& && 
$\{\lambda_i, \lambda_j, \lambda_k, \lambda_l\} \lambda_i \lambda_j
 d \lambda_k \wedge d \lambda_l$
\\
\\
\hline
\end{tabular}
\end{center}
\end{table}

\bibliographystyle{amsplain}
\bibliography{decomp}

\end{document}